\documentclass[USenglish]{article}

\usepackage[utf8]{inputenc}%(only for the pdftex engine)
\usepackage[small]{dgruyter}
\usepackage{microtype}
\usepackage{amssymb}
\usepackage{graphicx}
\usepackage{amsfonts,amssymb}
\usepackage{bm}
\usepackage{amsmath}
\usepackage{amsthm}
\usepackage{empheq}
\usepackage{cases}
\usepackage{arydshln}
\usepackage{dsfont}
\usepackage{subcaption}
\usepackage{algorithmicx}
\usepackage{algorithm}
\usepackage{algpseudocode}

\newtheorem{thm}{Theorem}
\newtheorem{exm}{Example}
\newtheorem{lem}{Lemma}
\newtheorem{rem}{Remark}
\newtheorem{pro}{Proposition}

\newtheorem{defn}{Definition}

\newcommand{\tr}{\mathrm{Tr}}

\newcommand{\blkdiag}{\mathrm{blkdiag}} 
\newcommand{\diag}{\mathrm{diag}} 
\newcommand{\im}{\mathrm{im}} 
\newcommand{\ibLambda}{{\it{\bar \Lambda }}}
\newcommand{\ihLambda}{{\it{ \hat \Lambda}}}

\begin{document}
\author*[1]{Xiaodong Cheng} 
\author[2]{Jacquelien M.A. Scherpen}
\author[3]{Harry L. Trentelman}
\runningauthor{X.~Cheng et al.}
\affil[1]{Eindhoven University of Technology}
\affil[2]{University of Groningen}
\affil[3]{University of Groningen}
\title{Reduced Order Modeling of Large-Scale Network Systems}
%\runningtitle{Introduction}
%  \subtitle{...}
%%%%%%%%%%%%%%%%%%%%%%%%%%%%%%%%%%%%%%%%%%%%%%%%%%%%%%%
\abstract{Large-scale network systems describe a wide class of complex dynamical systems composed of many interacting subsystems. A large number of subsystems and their high-dimensional dynamics often result in highly complex topology and dynamics, which pose challenges to network management and operation. This chapter provides an overview of reduced-order modeling techniques that are developed recently for simplifying complex dynamical networks. In the first part, clustering-based approaches are reviewed, which aim to reduce the network scale, i.e., find a simplified network with a fewer number of nodes. The second part presents structure-preserving methods based on generalized balanced truncation, which can reduce the dynamics of each subsystem. }
%%%%%%%%%%%%%%%%%%%%%%%%%%%%%%%%%%%%%%%%%%%%%%%%%%%%%%%
\classification[MSC~2010]{35B30, 37M99, 41A05, 65K99, 93A15, 93C05}
%\communicated{\today}
%  \dedication{...}
%\received{\today}
%\accepted{\today}
% \journalname{...}
% \journalyear{...}
% \journalvolume{..}
% \journalissue{..}
%\startpage{1}
% \aop
% \DOI{...}

\maketitle

\section{Introduction}
%%%%%%%%%%%%%%%%%%%%%%%%%%%%%%%%%%%%%%%%%%%%%%%%%%%%%%%%%%%%%%
Network systems, or multi-agent systems, are a class of structured systems composed of multiple interacting subsystems. 
In real life, systems taking the form of networks are ubiquitous, and the study of network systems has received compelling attention from many disciplines, see e.g., \cite{Newman2003Review,Newman2010NetworksIntroduction,mesbahi2010graph} for an overview. Coupled chemical oscillators, cellular and metabolic networks, interconnected physical systems, and electrical power grids
are only a few examples of such systems. 
To capture the behaviors and properties of network systems, graph theory is often useful \cite{godsil2013algebraic}. The interconnection structure among the subsystems can be represented by a \textit{graph}, in which \textit{vertices} and \textit{edges} represent the subsystems and the interactions among them, respectively.  
However, when network systems are becoming more large-scale, we have to deal with graphs of complex topology and nodal dynamics, which can cause great difficulty in transient analysis, failure detection, distributed controller design, and system simulation. 
From a practical point of view, it is always desirable to construct a reduced-order model to capture all the essential structure of the system avoiding too expensive computation. 
In the reduction of network systems, reduced-order models are designed not only to capture the main input-output feature of original complex network models but also to preserve the network structure such that they are usable for some potential applications, e.g., distributed controller design and sensor allocation.  

In the past few decades, a variety of theories and techniques of model reduction have been intensively investigated for generic dynamical systems. Techniques, including Krylov-subspace methods (also known as moment matching), balanced truncation, and Hankel norm approximation  \cite{Bai2002Krylov,Astolfi2010Moment,moore1981principal,glover1984all} provide us systematic procedures to generate reduced-order models that well approximate the input-output mapping of a high-dimensional system, see \cite{antoulas2005approximation,mor:vol1,mor:vol2} for an overview. However, when addressing the reduction of dynamical networks, the direct application of these methods may be not advisable, since they potentially destroy the network structure such that obtained reduced-order models could not have network feature any more. 
Structure-preserving model reduction is crucial for the application of network systems. Taking into account the two aspects of the complexity of network systems, namely, large-scale interconnection (i.e., a large number of subsystems) and high-dimensional subsystems, two types of problems are studied in the literature towards the approximation of network systems in a structure-preserving manner. 

The first problem aims to simplify the underlying network topology by reducing the number of nodes. The mainstream methods for this problem is based on \textit{graph clustering}, which is an unsupervised learning technique widely used in data science and computer graphics \cite{SurveyClustering,Schaeffer2007SurveyClustering}. For approximating dynamical networks, clustering-based methods basically follow a two-step process: the first is to partition the nodes into several nonoverlapping subsets (clusters), and then all the nodes in each cluster are aggregated into a single node. The aggregation step can be interpreted as a Petrov-Galerkin approximation using a clustering-based projection, see \cite{Schaft2014,Ishizaki2014,Monshizadeh2014}.
However, how to find the ``best'' clustering such that the approximation error is minimized still remains an open question. In \cite{Monshizadeh2014,jongsma2018model}, a particular clustering, called \textit{almost equitable partition} (AEP) is considered, which leads to an analytic $\mathcal{H}_2$ expression for the reduction error. However, finding AEPs itself is rather difficult and computationally expensive for general graphs. Clustering can also be found using the QR decomposition with column pivoting on the projection matrix obtained by the Krylov subspace method \cite{Petar2015CDC}. For undirected networks with tree topology, an asymptotically stable \textit{edge system} can be considered, which has a pair of diagonal generalized Gramian matrices for characterizing the importance of edges. Then, vertices linked by the less important edges are iteratively clustered \cite{Besselink2016Clustering}.  The notion of \textit{reducibility} is introduced in \cite{Ishizaki2014,ishizaki2015clustereddirected,Ishizaki2016dissipative} to characterize the uncontrollability of clusters. Using this notion, an upper bound for the network reduction error is established, which can determine the clustering. The works in \cite{ChengECC2016,ChengTAC2018MAS} extend the notion of dissimilarity for dynamical systems, where nodal behaviors are represented by the transfer functions mapping from external inputs to node states, and dissimilarity between two nodes are quantified by the norm of their behavior deviation. Then clustering algorithms, e.g., hierarchical clustering and K-means clustering, can be adapted to group nodes in such a way that nodes in the same cluster are more similar to each other than to those in other clusters  \cite{ChengTAC20172OROM,UmarCDC2019Clustering}. Subsequent research in \cite{ChengTAC20172OROM,ChengCDC20162OROM,ChengACOM2018Power,ChengTAC2019Digraph} show that the dissimilarity-based clustering method can also be extended to second-order networks, controlled power networks, and directed networks. In \cite{ChengECC2019weight,ChengTAC2020Weight}, a framework is presented on how to build a reduced-order model from a given clustering. The edge weights in the reduced network are parameterized so that an optimization problem is formulated to minimize the reduction error.

An alternative methodology to simplify the complexity of the network structure is based on time-scale analysis, and in particular, \textit{singular perturbation approximation}, see some of earlier works in \cite{Breusegem1991singular,phillips1981singular,biyik2008area}. Recently, this approach has also been extensively applied to biochemical systems and electric networks  \cite{Rao2013graph,anderson2011reduction,hancock2015simplified,ishizaki2013structured,Chow2013PowerReduction,Dorfler2013Kron,Monshizadeh2017Constant,romeres2013novel}. This class of approaches relies on the fact that the nodal states of network systems evolve over different time scales. Removing the vertices with fast states and reconnecting the remaining vertices with slow states will generate a reduced-order model that retains the low frequency behavior of the original network system. This methodology is closely related to the so-called \textit{Kron reduction} in electric networks \cite{Chow2013PowerReduction,Dorfler2013Kron,Monshizadeh2017Constant}, where the Schur complement of a graph Laplacian is taken that is again a Laplacian of a smaller-scale network. The singular perturbation approximation is capable of preserving the physical meaning of a network system. However, how to identify and separate fast/slow states is a crucial step in this approach, and its application is highly dependent on specific systems.

A network system can be simplified if the dimension of individual subsystems is reduced, which leads to the second research direction in reduced-order modeling of network systems see e.g., \cite{sandberg2009interconnected,Monshizadeh2013stability,ChengEJC2018Lure}. In this framework, the approximation is applied to each subsystem in
a way that certain properties of the overall network, such as synchronization and stability, are preserved. Relevant methods are developed using generalized Gramian matrices \cite{Dullerud2013CourseRobust} that allow for more freedom to preserve some desired structures than the standard Gramians. Networked nonlinear robustly synchronized Lur'e-type systems are reduced in \cite{ChengEJC2018Lure}, which shows that performing model reduction on the linear component of each nonlinear subsystems, allows for preserving the robust synchronization property of a Lur'e network. Techniques in \cite{Ishizaki2016dissipative,ChengAuto2019} can reduce the complexity of network structures and subsystem
dynamics \textit{simultaneously}. In \cite{Ishizaki2016dissipative}, the graph structure is simplified using clustering, while subsystems are reduced via some orthogonal projection. In contrast, \cite{ChengAuto2019} reduces graph structure and subsystem dynamics in a unified framework of generalized balanced truncation. Although a reduced-order system is obtained by balanced truncation does not necessarily preserve the network structure, a set of coordinates can be found in which the reduced-order model has a network interpretation. 

In this chapter, we will focus on the two problems of model reduction for linear network systems with diffusive couplings. In the aspect of simplifying network topology, we only review several clustering-based methods for space reasons. For the reduction of subsystems, we present the generalized balanced truncation as the main approach to perform a synchronization-preserving model reduction. The rest of this chapter is organized as follows. In Section~\ref{sec:Preliminaries}, we provide preliminaries on balanced truncation, semistable systems and necessary concepts in graph theory. The model of diffusively coupled networks is also introduced. In Section~\ref{sec:clustering}, we present clustering-based model reduction methods for simplifying network topology, and in Section \ref{sec:BTofNetwork}, the generalized balanced truncation approach is reviewed to reduce the dimension of subsystems. In Section \ref{sec:conclusion}, we glance at open problems and make some concluding remarks.

\textit{Notation:} The symbols $\mathbb{R}$ and $\mathbb{R}_+$ denote the set of real numbers and  real positive numbers, respectively. $I_n$ is the identity matrix of size $n$, and $\mathds{1}_n$ represents the vector in $\mathbb{R}^n$ of all ones. $\mathbf{e}_i$ is the $i$-th column of $I_n$. The cardinality of a finite set $\mathcal{S}$ is denoted by $|\mathcal{S}|$. $\tr(A)$, $\im(A)$, $\ker (A)$ denote the trace, image, and kernel of a matrix $A$, respectively. $\lVert G(s) \rVert_{\mathcal{H}_\infty}$ and $\lVert G(s) \rVert_{\mathcal{H}_2}$ represent the $\mathcal{H}_\infty$ norm and $\mathcal{H}_2$ norm of a transfer matrix $G(s)$.

\section{Preliminaries} \label{sec:Preliminaries}
In this section, we first briefly recapitulate the theory of balancing as a basis for the model reduction of linear control systems. New results for semistable systems and pseduo Gramians are also introduced. Moreover, we review some basic concepts from graph theory, which are then used for the modeling of network systems.

\subsection{Generalized Balanced Truncation} \label{sec:PreBT}
From \cite{Dullerud2013CourseRobust,antoulas2005approximation}, we review some basic facts on model reduction by using generalized balanced truncation. Consider a linear time-invariant system 
\begin{equation} \label{sys:ABC} 
\left\{
\begin{array}{l}
\dot{x} = A x  + B u, \\
y  = C x,  \\
\end{array}
\right.
\end{equation}
with $A \in \mathbb{R}^{n \times n}$, $B  \in \mathbb{R}^{n \times p}$, and $C \in \mathbb{R}^{q \times n}$, whose transfer function is given by
$G(s): = C (sI_n - A)^{-1} B$.
Let the system \eqref{sys:ABC} be asymptotically stable and minimal, i.e., $A$ is Hurwitz, the pair $(A, B)$ is controllable, and the pair $(C, A)$ is observable. Note that if a system \eqref{sys:ABC} is not minimal, we can always use the Kalman decomposition to remove the uncontrollable or unobservable states from the model \eqref{sys:ABC}. Thus, a minimal state-space realization can be obtained, of which the transfer function also is equal to $G(s)$.

For such a system \eqref{sys:ABC}, there always exits positive definite matrices, $P$ and $Q$  satisfying the following Lyapunov inequalities:
\begin{subequations}\label{eq:genGram}
	\begin{align}
	A P + P A^\top + B B^\top &\leq {0}, \\ 
	A^\top Q + Q A  + C^\top C &\leq {0}.
	\end{align}
\end{subequations}
Any $P$ and $Q$ as the solutions of \eqref{eq:genGram} are called \textit{generalized controllability and observability Gramians} of the system \eqref{sys:ABC} \cite{Dullerud2013CourseRobust}. When the equalities are achieved in \eqref{eq:genGram}, we obtain the standard controllability and observability Gramians, which become unique solutions of the Lyapunov equations \cite{antoulas2005approximation}.

Similar to the standard balancing, generalized balancing of the system \eqref{sys:ABC} amounts to finding a nonsingular matrix $T \in \mathbb{R}^{n\times n}$ such that $P$ and $Q$ are simultaneously diagonalized in the following way:
\begin{equation}
T P T^{\top} = T^{-\top} Q T^\top  = \Sigma :  = \text{diag}\left(\sigma_1, \sigma_2, \cdots, \sigma_n\right),
\end{equation}
where $\sigma_1 \geq \sigma_2 \geq \cdots \geq \sigma_n > 0$ are called   \textit{generalized Hankel singular values} (GHSVs) of system \eqref{sys:ABC}. Using $T$ as a coordinate transformation, we obtain a balanced realization of system \eqref{sys:ABC}, in which the state
components corresponding to the smaller GHSVs are relatively difficult to reach and observer and thus have less influence on the input-output behavior. 
Let the triplet $(\hat{A}, \hat{B}, \hat{C})$ be the $r$-dimensional reduced-order model (with $r \ll n$) obtained by truncating the states with the smallest GHSVs in the balanced system. Then, the reduced-order model $\hat{G}(s): = \hat{C} (sI_r - \hat{A})^{-1} \hat{B}$ preserves stability and moreover, an \textit{a priori} upper bound for the approximation error can be expressed in terms of the neglected GHSVs, i.e.,
\begin{equation} \label{eq:genbterr}
\lVert G(s) - \hat{G}(s) \rVert_{\mathcal{H}_\infty} \leq 2 \sum_{i=r+1}^{n} \sigma_i.
\end{equation} 

%%%%%%%%%%%%%%%%%%%%%%%%%%%%%%%%%%%%%%%
\subsection{Semistable Systems and Pseudo Gramians}
\label{sec:semistable}
Semistability is a more general concept than asymptotic stability as it allows for multiple zero poles in a system \cite{bhat1999lyapunov,hui2009semistability}. A linear system  $\dot{x} = Ax $ is \textit{semistable} if $\lim\limits_{t \rightarrow \infty} e^{At}$ exists for all initial conditions $x(0)$. The following lemma provides an equivalent condition for semistability. \\

\begin{lem} \cite{Semistable} \label{Ch6.lem:semisimplestable}
	A system $\dot{x} = Ax$ is semistable if and only if the zero eigenvalues of $A$ are {semisimple} (i.e., the geometric multiplicity of the zero eigenvalues coincides with the algebraic multiplicity), and all the other eigenvalues have negative real parts. 
\end{lem}

Let the triplet $(A, B, C)$ be a linear semistable system. The definition of semistablility implies that the transfer $G(s) = C (sI - A)^{-1}B$ is not necessarily in the $\mathcal{H}_2$ space, and the standard controllability and observability Gramians in \cite{antoulas2005approximation} are not well-defined in this case. Instead, we can define a pair of pseudo Gramians as follows \cite{ChengAuto2020Gramian}:   
\begin{equation}
\label{defn:PseudoGramians} 
\mathcal{P} = \int_{0}^{\infty} 
(e^{A t}-\mathcal{J}) BB^\top (e^{A^\top t}-\mathcal{J}^\top) \mathrm{d}t,
\
\mathcal{Q} = \int_{0}^{\infty} 
(e^{A^\top t}-\mathcal{J}^\top) C^\top C (e^{A t}-\mathcal{J}) \mathrm{d}t,
\end{equation}
where $\mathcal{J} : = \lim\limits_{t \rightarrow \infty} e^{At}$ is a constant matrix. The pseudo Gramians $\mathcal{P}$ and $\mathcal{Q}$ in \eqref{defn:PseudoGramians} are well-defined for semistable systems and can be viewed as a generalization of standard Gramian matrices for asymptotically stable systems. Furthermore, the pseudo Gramians can be computed as 
\begin{equation}
\mathcal{P} = \tilde{\mathcal{P}} - \mathcal{J} \tilde{\mathcal{P}} \mathcal{J}^\top, \quad \mathcal{Q} = \tilde{\mathcal{Q}} - \mathcal{J}^\top \tilde{\mathcal{Q}} {\mathcal{J}},
\end{equation}
where $\tilde{\mathcal{P}}$ and $\tilde{\mathcal{Q}}$ are arbitrary symmetric solution of the Lyapunov equations
\begin{align*}
A \tilde{\mathcal{P}} + \tilde{\mathcal{P}} A^\top + (I-\mathcal{J})BB^\top(I-\mathcal{J}^\top) = 0,  \\
A^\top \tilde{\mathcal{Q}} + \tilde{\mathcal{Q}} A  + (I-\mathcal{J}^\top)C^\top C(I-\mathcal{J}) = 0,
\end{align*}
respectively. The pseudo Gramians lead to a characterization of the $\mathcal{H}_2$ norm of a semistable system.   
\begin{thm} \cite{ChengAuto2020Gramian} \label{thm:H2semi}
	Consider a semistable system with the triplet $(A, B, C)$. Then, $G(s): = C(sI-A)^{-1}B \in \mathcal{H}_2$  if and only if $C \mathcal{J} B = 0$. Furthermore, if $\| G(s) \|_{\mathcal{H}_2}$ is well-defined, then
	\begin{equation}
	\lVert G(s) \rVert_{\mathcal{H}_2}^2 = \tr (C \mathcal{P} C^\top) = \tr (B^\top \mathcal{Q} B).
	\end{equation}
\end{thm}

%%%%%%%%%%%%%%%%%%%%%%%%%%%%%%%%%%%
\subsection{Graph Theory}

The concepts in graph theory are instrumental in analyzing network systems \cite{mesbahi2010graph}. The interconnection structure of a network is often characterized by a graph $\mathcal{G}$ that consists of a finite and nonempty node set  $\mathcal{V}: = \{1, 2, ... , n\}$ and an edge set $\mathcal{E} \subseteq \mathcal{V} \times \mathcal{V}$.
Each element in $\mathcal{E}$ is an ordered pair of elements of $\mathcal{V}$, and we say that the edge is directed from vertex $i$ to vertex $j$ if $(i,j) \in \mathcal{E}$. This leads to the definition of the \textit{incidence matrix} $R \in \mathbb{R}^{n \times |\mathcal{E}|}$:
\begin{equation} \label{eq:incidence}
[R]_{ij} = \begin{cases} 
+1 & \text{if edge}~(i,j)\in \mathcal{E}, \\ 
-1 & \text{if edge}~(j,i)\in \mathcal{E}, \\ 
0 & \text{otherwise}. \end{cases} 
\end{equation}
If each edge is assigned a positive value (weight), the graph $\mathcal{G}$ is \textit{weighted}, and a \textit{weighted adjacency matrix} $\mathcal{W}$ can be defined such that $w_{ij} = [\mathcal{W}]_{ij}$ is positive if there exists a directed edge from node $j$ to node $i$, i.e., $(j,i) \in \mathcal{E}$, and $w_{ij}=0$ otherwise. A (directed) graph $\mathcal{G}$ is called \textit{undirected} if $\mathcal{W}$ is symmetric. 
An undirected graph $\mathcal{G}$ is called \textit{simple}, if $\mathcal{G}$ does not contain self-loops (i.e., $\mathcal{E}$ does not contain edges of the form $(i,i)$, $\forall~i$), and there exists only one undirected edge between any two distinct nodes. Two distinct vertices $i$ and $j$ are said to be \textit{neighbors} if there exists an edge between $i$ and $j$, and the set $\mathcal{N}_i$ denotes all the neighbors of node $i$.

The \textit{Laplacian matrix} $L \in \mathbb{R}^{n \times n}$ of a weighted graph $\mathcal{G}$ is defined as
\begin{equation} \label{defn:Laplacian}
[L]_{ij} = \left\{ \begin{array}{ll} 
\sum_{j \in \mathcal{N}_i}^{} w_{ij}, & i = j\\
-w_{ij}, & \text{otherwise.}
\end{array}
\right.
\end{equation}
Furthermore, we can define an undirected graph Laplacian using an alternative formula:
\begin{equation} \label{eq:LRWR}
L = R W R^\top,
\end{equation}
where $R$ is an incidence matrix obtained by assigning an arbitrary orientation to each edge of $\mathcal{G}$, and $W := \diag(w_1, w_2, ..., w_{|\mathcal{E}|})$ with $w_k$ the weight associated with the edge $k$, for each $k = 1,2,..., |\mathcal{E}|$.
\vspace*{12pt}

\begin{rem} \label{rem:Lpro}
	If $\mathcal{G}$ is a simple undirected connected graph, the associated Laplacian matrix $L$ has the following structural properties:
	\begin{enumerate}
		\item ${L}^\top = {L}  \geq 0$ 
		\item $\ker(L)  = \im(\mathds{1})$;
		\item $L_{ij} \leq 0$ if $i \ne j$, and $L_{ij} > 0$ otherwise.
	\end{enumerate}
	Conversely, any real square matrix satisfying the above conditions can be interpreted as a Laplacian matrix that uniquely represents a simple undirected connected graph.
\end{rem}

\subsection{Network Systems}

In this chapter, we mainly focus on an important class of networks, namely \textit{consensus networks}, where subsystems are interconnected via \textit{diffusive couplings}. Various applications, including formation control of mobile vehicles, synchronization in power networks, and balancing in chemical kinetics, involve the concept of consensus networks \cite{ren2005survey,consensus2010li,fax2001graph,dorfler2014synchronization,CucuzzellaTCST2018Consensus,Trip2018distributed,Rao2015balancing}.

Here, we consider a network system in which the interconnection
structure is represented by a simple weighted undirected  graph with the node set $\mathcal{V}= \{1, 2, ..., n\}$. The dynamics of each vertex (agent) is described by 
\begin{equation} \label{sysagent}
{\bm{\Sigma}_i}:\left \{
\begin{array}{l}
\dot{x}_i = A x_i + B v_i, \\
\eta_i = C x_i,
\end{array}
\right.    
\end{equation}
where $x_i \in \mathbb{R}^{\ell}$, $v_i \in \mathbb{R}^{{m}}$ and $\eta_i \in \mathbb{R}^{{m}}$ are the state, control input and output of node $i$, respectively. The $n$ subsystems are interconnected such that
\begin{equation} \label{eq:coupling}
m_i {v}_i = -\sum_{j \in \mathcal{N}_i}^{ } w_{ij} \left(\eta_i-\eta_j\right)
+ \sum_{j=1}^{p} f_{ij} u_j,
\end{equation}
where $m_i \in \mathbb{R}_{+}$ denotes the weight of node $i$.
In \eqref{eq:coupling}, the first term on the left is referred to as \textit{diffusive coupling}, where $w_{ij} \in \mathbb{R}_{+}$ is the entry of the adjacency
matrix $[\mathcal{W}]_{ij}$ standing for the intensity of the coupling between nodes $i$ and $j$. The second term indicates the influence from external input
$u_j$, where the value of $f_{ij} \in \mathbb{R}$ represents the amplification of $u_j$ acting on vertex $i$.  Let $F\in \mathbb{R}^{n \times p}$ be a matrix with $[F]_{ij} = f_{ij}$, and we introduce the external outputs as $y_i = \sum_{j=1}^{n} [H]_{ij} \eta_j$, with $y_i \in \mathbb{R}^{m}$, as the $i$-th external output of the network. 
We then represent the network system in compact form as
\begin{equation} \label{sysh}
\bm{\Sigma}:  
\left \{
\begin{array}{rcl}
(M \otimes I)\dot{x} &=& \left(M \otimes A - L \otimes BC\right)x + (F \otimes B)u 
\\
y&=&(H \otimes C)x.
\end{array}
\right.
\end{equation}
with joint state vector $x^\top := \left[x_1^\top\ x_2^\top\ ...\ x_n^\top \right] \in \mathbb{R}^{n \ell}$, external control input $u^\top := \left[{u}_1^\top\ u_2^\top\ ... \ u_p^\top \right]\in \mathbb{R}^{pm}$ and external output $y = \left[{y}_1^\top\ y_2^\top\ ... \ y_q^\top \right] \in \mathbb{R}^{qm}$. 
$M : = \diag(m_1, m_2, ..., m_n) > 0$, and
$L \in \mathbb{R}^{n \times n}$ is the graph Laplacian matrix that characterizes the coupling structure among the subsystems. In many studies of undirected networks, the matrix $M = I_n$ is considered.

The simplest scenario in network systems is that all
the vertices are just single-integrators, i.e., $m_i \dot{x}_i = v_i$ with $x_i \in \mathbb{R}$. Then, the model of a networked single-integrator system can be formed by taking $A=0$ and $B=C=1$ in \eqref{sysh}, which leads to 
\begin{equation} \label{sys}
\left \{ \begin{array}{rcl}
M\dot{x} &=& -Lx + Fu, \\
y &=& Hx.
\end{array} \right.
\end{equation}
A variety of physical systems are of this form, such as {mass-damper systems} and {single-species reaction networks} \cite{Schaft2014}. Note that the system \eqref{sys} is call \textit{semistable} \cite{bhat1999lyapunov}, since $L$ has a simple zero eigenvalue.

%Clustering-based model reduction methods for the system \eqref{sys} can be found in e.g., \cite{Schaft2014,Monshizadeh2014,ChengECC2016}.   

An important issue in the context of diffusively coupled networks is \textit{synchronization}. The system $\bm{\Sigma}$ in \eqref{sysh} achieves synchronization if, for any initial conditions, the zero input response of \eqref{sysh} satisfies
\begin{equation}
\lim\limits_{t \rightarrow \infty}  \left[x_i(t) - x_j(t)\right]=0, \  \text{for all} \  i,j \in \mathcal{V}.
\end{equation}
Using the property of $L$ in Remark~\ref{rem:Lpro}, it is clear that the single-integrator network in \eqref{sys} can reach synchronization. However, for the general form of \eqref{sysh}, we need take into account the subsystems as well.
Denote by $0 = \lambda_1 < \lambda_2 \leq \cdots \leq \lambda_{n}$  the eigenvalues of the matrix $M^{-1}L$ in ascending order. A sufficient and necessary condition for the synchronization of a network consisting of agents as in \eqref{sysagent} is found in e.g. \cite{consensus2010li}.
\vspace*{12pt}

\begin{lem}  
	\label{lem:ConsCond}
	The multi-agent system $\bm{\Sigma}$ in \eqref{sysh} achieves {synchronization} if and only if $A-\lambda_{k} BC$ is Hurwitz, for all $k \in \{2,3,...,n\}$.  
\end{lem}

\section{Clustering-Based Model Reduction}
\label{sec:clustering}
%%%%%%%%%%%%%%%%%%%%%%%%%%%%%%%%%%%%%%%%%%%%%%%%%%%%%%%%%%%%%%

In this section, we introduce clustering-based methods that combine the Petrov-Galerkin projection with graph clustering. A reduced-order network model can be constructed by using the characteristic matrix of a graph clustering. Moreover, we will also briefly recap some other clustering-based methods.

Graph clustering is an important notion in graph theory  \cite{godsil2013algebraic}. Consider a connected undirected graph $\mathcal{G} = \left( \mathcal{V}, \mathcal{E}\right)$. A \textit{graph clustering} of $\mathcal{G}$ is to divide its vertex set $\mathcal{V}$ (with $|\mathcal{V}| = n$) into $r$ nonempty and disjoint subsets, denoted by  $   \mathcal{C}_1,\mathcal{C}_2,...,\mathcal{C}_r$, where $\mathcal{C}_i$ is called a \textit{cluster} (or a \textit{cell} of $\mathcal{G}$). 
\begin{defn}    
	The {characteristic matrix} of the clustering $\{  \mathcal{C}_1,\mathcal{C}_2,...,\mathcal{C}_r\}$ is characterized by the binary matrix $\Pi \in \mathbb{R}^{n \times r}$ as
	\begin{equation} \label{partition}
	[\Pi]_{ij} := \left\{ \begin{array}{ll}
	1, & \text{if vertex $i \in \mathcal{C}_j$,}\\ 0, & \text{otherwise.}
	\end{array}
	\right.
	\end{equation}
\end{defn}
\vspace*{12pt}

Note that each vertex is assigned to a unique cluster. Therefore, each row of the characteristic matrix $\Pi$  has exactly one nonzero element, and the number of nonzero elements in each column of $\Pi$ is the number of vertices in the corresponding cluster. Specifically, we have 
\begin{equation} \label{eq:Piproperty}
\Pi \mathds{1}_r = \mathds{1}_n~\text{and}~\mathds{1}^\top_n \Pi=\left[|\mathcal{C}_1|, |\mathcal{C}_2|,...,|\mathcal{C}_r| \right].
\end{equation}
It is worth noting that for any given undirected graph Laplacian $L$, the reduced matrix $\Pi^\top L \Pi$ is also a Laplacian matrix, representing an undirected graph of smaller size. This important property allows for a structure-preserving model reduction of network systems using $\Pi$ for the Petrov-Galerkin projection.

Let $\bm{\Sigma}$ in \eqref{sysh} be a network system with underlying graph $\mathcal{G}$ of $n$ vertices. To formulate a reduced-order network model of $r$ dimension, we first find a graph clustering that partitions the vertices of $\mathcal{G}$ into $r$ clusters.  Then we use the characteristic matrix of the clustering as a basis that projects the state-space of $\bm{\Sigma}$ to a reduced subspace. Specifically,
a reduced-order model of $\bm{\Sigma}$ can be constructed via the Petrov-Galerkin projection framework as 
\begin{equation} \label{sysrh}
\bm{\hat{\Sigma}} : \left \{
\begin{array}{rcl}
(\hat{M} \otimes I)\dot{z} & = & \left(\hat{M} \otimes A - \hat{L} \otimes B\right)z + (\hat{F} \otimes B) u, \\
\hat{y} & = &  (\hat{H} \otimes C) z,
\end{array}
\right.
\end{equation}
where $\hat{M} := \Pi^\top M \Pi \in \mathbb{R}^{r \times r}$, $\hat{L} := \Pi^\top {L} \Pi \in \mathbb{R}^{r \times r}$, $\hat{F}=\Pi^\top F$ and $\hat{H} = H\Pi$. The new state vector $z^\top: =  \left[z_1^\top\ z_2^\top\ ...\ z_r^\top \right] \in \mathbb{R}^{r \ell}$, where each component $z_i \in \mathbb{R}^{\ell}$ represents an estimated dynamics of all the vertices in the $i$-th cluster, and $\hat{x} = (\Pi \otimes I) z \in \mathbb{R}^{n \ell}$ can be an approximation of the original state $x$. For the single-integrator network system \eqref{sys}, clustering-based projection yields the reduced-order model as
\begin{equation} \label{sysr}
\left \{
\begin{array}{rcl}
\hat{M} \dot{z} &=& - \hat{L} z + \hat{F} u, \\
\hat{y} &=& \hat{H} z.
\end{array}
\right.
\end{equation} 
In the reduced-order models in \eqref{sysrh} and \eqref{sysr}, $\hat{M}$ is a positive diagonal matrix, and $\hat{L} \in \mathbb{R}^{r \times r}$ is a Laplacian matrix representing a graph of a fewer number of vertices. More preciously, $\hat{M}$ and $\hat{L}$ can be computed as
\begin{equation}
[\hat{M}]_{kk} = \sum_{i \in \mathcal{C}_k} m_i, 
\quad [\hat{L}]_{kl} = 
\begin{cases}
\sum_{i \in \mathcal{C}_k, j \in \mathcal{C}_l}^{} [L]_{ij}, & k \ne l; 
\\
\sum_{i \in \mathcal{C}_k}^{} [L]_{ij}, & k = l.
\end{cases}
\end{equation}
Clearly, the reduced-order models in \eqref{sysrh} and \eqref{sysr} preserve the network structure and thus can be interpreted as simplified dynamical networks with diffusive couplings. The following example then illustrates the physical meaning of a projected reduced-order network model.
\\

\begin{exm} \label{exm1}
	Consider a mass-damper system in Fig.~\ref{example1}, left inset, where ${u}_1$, $u_2$ are external forces acting on the 1st and 4th mass blocks. Suppose that all the masses are identical, then we model the network system in the form of \eqref{sys} with   
	\begin{equation*} 
	M = I_5, \
	L=\left[                
	\begin{array}{rrrrr}
	6  & -3&  0 & -2 & -1\\
	-3 &  4&  -1&  0 & 0\\
	0  &  -1&  6&  -2 & -3\\
	-2 &  0&  -2&  5 & -1\\
	-1 &  0&  -3&  -1 & 5\\
	\end{array}
	\right], \   
	%%%%%%%%%%%%%%%%%%%%%%%%%%
	F= \left[\begin{matrix}
	1 & 0 \\   
	0 & 0 \\
	0 & 0 \\
	0 & 1 \\
	0 & 0 \\
	\end{matrix}\right].            
	\end{equation*}
	The off-diagonal entry $-[L]_{ij}$ represents the damping coefficient of the edge $(i,j)$.  Let $\{\mathcal{C}_1, \mathcal{C}_2, \mathcal{C}_3\} = \left\{ \{1,2\}, \{3,5\},\{4\}\right\}$ be the clustering of the network, which leads to
	\begin{equation*}
	\Pi=\left[                 
	\begin{matrix} 
	1 & 1 &  0 & 0 & 0\\
	0 & 0 &  1 & 0 & 1\\
	0 & 0 &  0 & 1 & 0\\
	\end{matrix}
	\right]^\top.   
	\end{equation*}
	A reduced-order network model is obtained as in \eqref{sysr} with
	\begin{equation*}  
	\begin{split}
	\hat{M}= \left[                
	\begin{array}{ccc}   
	2 & 0 & 0\\
	0 & 2 & 0\\
	0  & 0 & 1\\
	\end{array}
	\right],\ 
	&\hat{L}  = \left[                
	\begin{array}{rrr}   
	4 &  -2  &   -2\\
	-2  &   5   & -3\\
	-2 &   -3  & 5
	\end{array}
	\right], \ 
	%%%%%%%%%%%%%%%%%%%%%%%%%%%%%%%
	\hat{F} =\left[\begin{matrix}
	1 & 0  \\
	0 & 1 \\
	0 & 1
	\end{matrix}\right].
	\end{split} 
	\end{equation*}
	It is clear that each mass in the reduced network is
	equal to the sums of the masses in the corresponding cluster. Moreover, the structure of a Laplacian matrix is retained, which allows for a physical interpretation of the reduced model, as shown in Fig. \ref{example1}, right inset.    
	\begin{figure}[!tp]\centering
		\includegraphics[scale=.6]{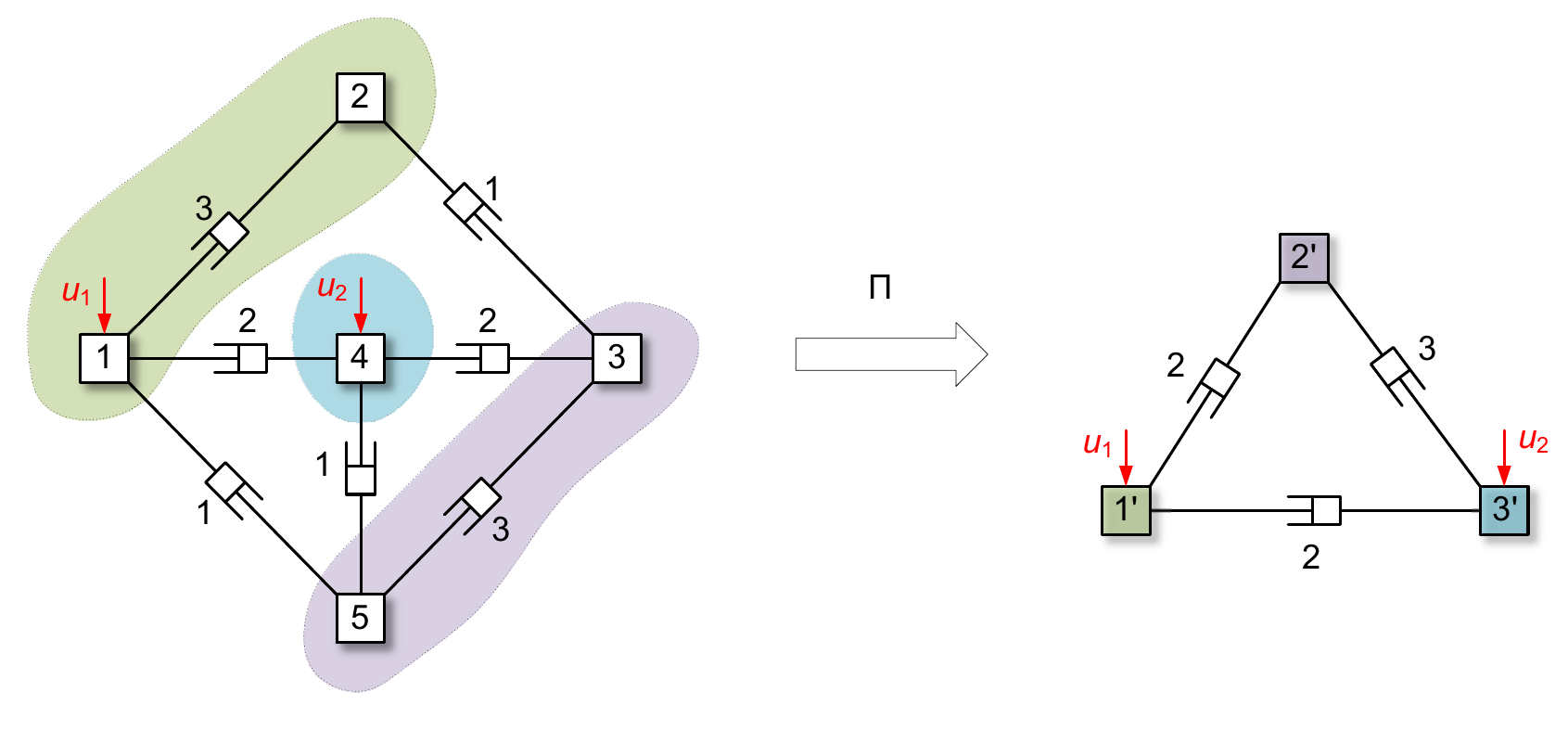}     
		\caption{An illustrative example of clustering-based model reduction for a mass-damper network system.}
		\label{example1}
	\end{figure}
\end{exm}

Next, the property of the reduced-order models in \eqref{sysrh} and \eqref{sysr} are discussed. First, it is clear that system \eqref{sysr} preserves the synchronization property. Moreover, the following result holds.
\begin{lem} \cite{ChengECC2016,ChengTAC2018MAS}
	Consider the single-integrator networks in \eqref{sys} and \eqref{sysr}. The impulse responses of the two systems satisfy
	\begin{equation}
	\lim\limits_{t \rightarrow \infty} y(t) 
	= \lim\limits_{t \rightarrow \infty} \hat{y}(t) 
	= \frac{H \mathds{1}_n \mathds{1}_n^\top F}{\mathds{1}_n^\top M \mathds{1}_n}
	\end{equation}
\end{lem}
Denote 
\begin{equation} \label{eq:TFS}
S: = H (sI_n + L)^{-1} F, \quad \hat{S}: = \hat{H} (sI_r + \hat{L})^{-1} \hat{F}
\end{equation}
This lemma implies the reduction error $\| S - \hat{S} \|_{\mathcal{H}_2}$ is well defined, for any clustering $\Pi$.
For the reduced-order network system \eqref{sysrh}, the analysis of synchronization and reduction error is more complicated, since the subsystem dynamics will also be involved. 
Denote
\begin{subequations} \label{eq:TF}
	\begin{align}
	G(s): &= (H \otimes C) \left[ M \otimes (sI_\ell - A) + L \otimes BC \right](F \otimes B),  \\
	\hat{G}(s): &= (\hat{H} \otimes C) \left[ \hat{M} \otimes (sI_\ell - A) + \hat{L} \otimes BC \right](\hat{F} \otimes B),
	\end{align}
\end{subequations} 
as the transfer matrices of the systems \eqref{sysh} and \eqref{sysrh}, respectively. Generally, $G(s) - \hat{G}(s)$ is not guaranteed to be stable.
However, a theoretical guarantee can be obtained if the subsystem $(A, B, C)$ in \eqref{sysagent} are \textit{passive} \cite{JanWillems1976}, namely, there exists a positive definite $K$ such that     
\begin{equation} \label{eq:passive}
A^\top K + K A \leq 0, \ \text{and} \ C^\top = B K.
\end{equation}
In this case, we have synchronization property and bounded reduction error for the system \eqref{sysrh}. 
\begin{thm}\label{thm:syshsyn}
	Consider the subsystem $(A, B, C)$ in \eqref{sysagent}, which is passive and minimal. Then the following statements hold.
	\begin{enumerate} 
		\item The original network system \eqref{sysh} achieves synchronization for any $L$ representing an undirected connected graph \cite{scardovi2008synchronization,chopra2012output}.
		
		\item The reduced-order network system \eqref{sysrh} achieves synchronization for any clustering $\Pi$ \cite{Besselink2016Clustering,ChengTAC2018MAS}. 
		
		\item $G(s) - \hat{G}(s) \in \mathcal{H}_2$, for any any clustering $\Pi$ \cite{Besselink2016Clustering,ChengTAC2018MAS}. 
	\end{enumerate}
\end{thm}

In the framework of clustering-based projection, the approximation error $\| G(s) - \hat{G}(s) \|_{\mathcal{H}_2}$  only depends on the choice of graph clustering. Thus, it is a crucial problem in this framework to select a suitable clustering such that the obtained reduced-order model \eqref{sysrh} can well approximate the behavior of the original network system \eqref{sysh}. In the following subsections, we review several cluster selection methods. 

\subsection{Almost Equitable Partitions} \label{XD.sec:aep}

It is suggested in \cite{Monshizadeh2014} to place those vertices that connect to the rest of the network in a similar fashion into the same cluster. This idea leads to
a special class of graph clusterings, namely \textit{almost equitable partitions}. 
\begin{defn}
	Let $\mathcal{G} = (\mathcal{V}, \mathcal{E})$ be a weighted undirected graph. A graph clustering $\{\mathcal{C}_1, \mathcal{C}_2, ..., \mathcal{C}_r\}$ is called an almost equitable partition if for each $\mu, \nu \in \{1,2,...,r\}$ with $\mu \ne \nu$, it holds that
	$
	\sum_{k \in \mathcal{C}_\nu}^{}w_{ik} = \sum_{k \in \mathcal{C}_\nu}^{}w_{jk}
	$, $\forall~i,j \in \mathcal{C}_\mu$, where $w_{ij}$ denotes the $(i,j)$-the entry of the adjacency matrix of $\mathcal{G}$.
\end{defn}

If $\{\mathcal{C}_1, \mathcal{C}_2, ..., \mathcal{C}_r\}$ is an almost equitable partition, its characteristic matrix $\Pi$ has the key property that $\im(\Pi)$ is $L$-invariant \cite{Monshizadeh2014}, i.e., 
$
L~\im(\Pi) \subseteq \im(\Pi).
$

Consider the single-integrator network in \eqref{sys} with $\mathcal{V}$ the vertex set. In the context of leader-follower networks, a subset of vertices $\mathcal{V}_L \subseteq \mathcal{V}$ are the leaders, with $|\mathcal{V}_L| = p$, which are controlled by external inputs. Moreover, $F \in \mathbb{R}^{n \times p}$ in \eqref{sys} is the binary matrix such that $[F]_{ij} = 1$ if vertex $i$ is the $j$-th leader, and $[F]_{ij} = 0$ otherwise.
Assume that the output of \eqref{sys} is given as
\begin{equation} \label{eq:outputAEP}
y =  H x = W^{\frac{1}{2}} R^\top x,
\end{equation}
where $R$ is the incidence matrix of $\mathcal{G}$, and $W$ is the edge weight matrix defined in \eqref{eq:LRWR}.  Then, the output of the reduced network model  \eqref{sysr} is
obtained as $\hat{y} = \hat{H}x =  W^{\frac{1}{2}} R^\top \Pi x$ with $\Pi$ the characteristic matrix of the given almost equitable partition. Using the property of the output matrices: $H^\top H = L$ and $\hat{H}^\top \hat{H} = \hat{L}$,  an explicit
$\mathcal{H}_2$ error can be derived, which is characterized by the cardinalities of the clusters
containing leaders \cite{Monshizadeh2014,jongsma2018model}.

\begin{thm}
	Consider the network system  \eqref{sys} with output defined in \eqref{eq:outputAEP}.
	Let $\Pi$ be the characteristic matrix of an almost equitable partition of the underlying graph: $\{\mathcal{C}_1, \mathcal{C}_2, ..., \mathcal{C}_r\}$. Denote $S$ and $\hat{S}$ as the transfer matrices of \eqref{sys} and \eqref{sysr}, respectively.
	Then, we have 
	\begin{equation}
	\dfrac{\| S -  \hat{S} \|_{\mathcal{H}_2}^2}{\| S\|_{\mathcal{H}_2}^2}
	=
	\dfrac{\| S   \|_{\mathcal{H}_2}^2- \|   \hat{S} \|_{\mathcal{H}_2}^2}{\| S\|_{\mathcal{H}_2}^2}
	= 
	\dfrac{\sum_{i=1}^{p}\left(1 - \frac{1}{|C_{k_i}|}\right)}{p\left(1 - \frac{1}{n}\right)},
	\end{equation}
	where 
	$n = |\mathcal{V}|$, $p = |\mathcal{V}_L|$, and 
	$k_i$ is the integer index such that the $i$-th leader is within $\mathcal{C}_{k_i}$.
\end{thm}

In \cite{jongsma2018model},  a formula for the $\mathcal{H}_\infty$-error is also derived by assuming a specific output $y = L x$ in \eqref{sys}. If the network \eqref{sys} is clustered according to an almost equitable partition $\{\mathcal{C}_1, \mathcal{C}_2, ..., \mathcal{C}_r\}$, then we have
\begin{equation*}
\| S - \hat{S} \|^2_{\mathcal{H}_\infty} = 
\begin{cases}
\max\limits_{1 \leq i \leq q} \left(1 - \frac{1}{|\mathcal{C}_{k_i}|}\right) & \text{if the leaders are in different clusters},\\
1 & \text{otherwise},
\end{cases}
\end{equation*} 
with $k_i$ the integer such that the $i$-th leader is within $\mathcal{C}_{k_i}$.

More results on model reduction methods based on almost equitable partitions can be found in \cite{jongsma2018model}, where network systems of the form \eqref{sysh} with symmetric subsystems are also discussed.
\vspace{20pt}

\subsection{Clustering of Tree Networks} \label{XD.sec:tree}

If the underlying graph of the considered network model \eqref{sysh} is a tree, we can resort to the model reduction procedure proposed in \cite{Besselink2016Clustering}.
%According to \cite{JanWillems1976}, the passivity of the subsystems \eqref{sysagent} means that there exists a positive definite $K$ such that 
%\begin{equation} \label{eq:passive}
%    A^\top K + K A \leq 0, \ \text{and} \ C^\top = B K.
%\end{equation}
Consider the network model $\bm{\Sigma}$ in \eqref{sysh}, where the subsystems are passive and minimal, and the Laplacian matrix $L$ represents an undirected tree graph $\mathcal{T}$. 
Note that if $\mathcal{T}$ contains $n$ vertices, then it has $n-1$ edges. Relevant to \eqref{eq:LRWR}, an \textit{edge Laplacian} is defined:
\begin{equation}
L_\mathrm{e} = R^\top R W,
\end{equation}
where $R \in \mathbb{R}^{n \times (n -1)}$ is the
oriented incidence matrix of $\mathcal{T}$, and $W \in \mathbb{R}^{(n-1) \times (n-1)}$ is the edge weight matrix. It is not hard to see that $L_\mathrm{e}$ has all eigenvalues real and positive, and these eigenvalues coincide to the nonzero eigenvalues of $L$. 

Let $M = I_n$ in \eqref{sysh}, and an edge system can be defined as
\begin{equation} \label{sys:edge}
\bm{\Sigma}_\mathrm{e}: \begin{cases}
\dot{x}_\mathrm{e} = (I_{n-1} \otimes A - L_\mathrm{e} \otimes BC) x_\mathrm{e} + (F_\mathrm{e} \otimes B) u, \\
y_\mathrm{e} = (H_\mathrm{e} \otimes C) x_\mathrm{e},
\end{cases}
\end{equation}
where $x_\mathrm{e} = (R^\top \otimes I) x \in \mathbb{R}^{(n-1)\ell}$, $F_\mathrm{e} = R^\top F$, and $H_\mathrm{e} = H R W L_\mathrm{e}^{-1} $. 
A dual edge system is also introduced with a different realization as
\begin{equation} \label{sys:edgedual}
\bm{\Sigma}_\mathrm{f}: \begin{cases}
\dot{x}_\mathrm{f} = (I_{n-1} \otimes A - L_\mathrm{e} \otimes BC) x_\mathrm{f} + (F_\mathrm{f} \otimes B) u, \\
y_\mathrm{e} = (H_\mathrm{f} \otimes C) x_\mathrm{f},
\end{cases}
\end{equation}
with $x_\mathrm{f} = (L_\mathrm{e}^{-1} \otimes I) x_\mathrm{e}$, $F_\mathrm{f} = L_\mathrm{e}^{-1} F_\mathrm{e}$, and $H_\mathrm{f} = H R W$. 

Assuming that $(A, B, C)$ is passive and minimal, the system \eqref{sysh} achieves synchronization from Theorem~\ref{thm:syshsyn}, which means that $A- \lambda_{k} BC$ is Hurwitz for all nonzero eigenvalues $\lambda_{k}$ of graph Laplacian matrix $L$. This further implies that both $\bm{\Sigma}_\mathrm{e}$ and $\bm{\Sigma}_\mathrm{f}$ are asymptotically stable.  As a result, generalized controllability and observability Gramians of the edge systems \eqref{sys:edge} and \eqref{sys:edgedual} can be analyzed. 
\\

\begin{lem} \cite{Besselink2016Clustering}\label{lem:XYLe}
	Consider the edge systems \eqref{sys:edge} and \eqref{sys:edgedual} of a tree network. There exist matrices $X >0$ and $Y >0$ such that the following inequalities hold:
	\begin{eqnarray}
	-L_\mathrm{e} X - X L_\mathrm{e}^\top + R^\top FF^\top R \leq 0, \\
	-L_\mathrm{e}^\top Y - Y L_\mathrm{e} + WR^\top FF^\top RW \leq 0. 
	\end{eqnarray}
	Moreover, $P_\mathrm{e} : = X \otimes K^{-1}$ and $Q_\mathrm{f} : = Y \otimes K$ are a generalized controllability Gramian of $\bm{\Sigma}_\mathrm{e}$ in \eqref{sys:edge} and a generalized observability Gramian of $\bm{\Sigma}_\mathrm{f}$ in \eqref{sys:edgedual}, respectively, where $K$ satisfies \eqref{eq:passive} for the passive subsystems. 
\end{lem}

According to \cite{Besselink2016Clustering}, the matrices $X$ and $Y$ can be chosen to admit a diagonal structure: 
\begin{equation}\label{eq:XY}
X = \diag (\xi_1, \xi_2, ..., \xi_{n-1}), \quad Y = \diag (\eta_1, \eta_2, ..., \eta_{n-1}),
\end{equation}
where the ordering $\xi_i \eta_i \geq \xi_{i+1} \eta_{i+1}$ is imposed. Note that $X$ and $Y$ imply the controllability and observability properties of the edges, respectively, the value of $\xi_i \eta_i$ can be viewed as an indication for the importance of the $i$-th edge. Following a similar reasoning as balanced truncation in Section~\ref{sec:PreBT}, removing the edges according to the value of $\xi_i \eta_i$ is meaningful. In \cite{Besselink2016Clustering}, a graph clustering procedure is presented to recursively aggregate the two vertices connected by the least important edge. Furthermore, an \textit{a priori} upper bound on the approximation error in terms of the $\mathcal{H}_\infty$ norm can be derived.

\begin{thm}
	Consider the networked system in \eqref{sys} with $M = I_n$. Assume each subsystem is minimal and passive, and the underlying graph is a tree. Let \eqref{sysrh} be the $r$-th order reduced network system obtained by aggregating the vertices connecting by the least important edges of the original network. Then, the following error bound holds.
	\begin{equation}
	\| G(s) - \hat{G}(s) \|_{\mathcal{H}_\infty} \leq 
	2 \left( \sum_{i = r }^{n-1} [L_\mathrm{e}^{-1}]_{ii} \sqrt{\xi_i \eta_i}\right),
	\end{equation}
	where $G(s)$ and $\hat{G}(s)$ are transfer matrices in \eqref{eq:TF}, and  $[L_\mathrm{e}^{-1}]_{ii}$ denotes the $i$-th diagonal entry of the matrix $L_\mathrm{e}^{-1}$, and $\xi_i$ and $\eta_i$ are the diagonal entries of $X$ and $Y$ in \eqref{eq:XY}, respectively. 
\end{thm}

Note that the proposed method in \cite{Besselink2016Clustering} heavily relies on the assumption of tree topology. For networks with general topology, applying this method would be challenging, since there may not exist edge systems as in \eqref{sys:edge} and \eqref{sys:edgedual}, which admit diagonal Gramians as in \eqref{eq:XY}.  

\subsection{Dissimilarity-Based Clustering} \label{XD.sec:projection}

The methods in Section \ref{XD.sec:aep} and Section \ref{XD.sec:tree} rely either on a special graph clustering or on a specific topology. In this section, we review a  dissimilarity-based method, which can be performed to reduce more general network systems. Clustering of data points in data science is usually based
on some similarity measure in terms of vector norms. 
To cluster a dynamical network, we can extend the concept of \textit{dissimilarity} using the function norms, which serves as a metric for quantifying how differently two distinct vertices (subsystems) behave \cite{ChengECC2016,ChengTAC2018MAS}. 
\begin{defn} \label{defn:dissim}
	Consider a network system in \eqref{sysh} or \eqref{sys}, the {dissimilarity} between vertices $i$ and $j$ is defined as 
	\begin{equation} \label{eq:dissim}
	\mathcal{D}_{ij} := \lVert \eta_i(s) - \eta_j(s) \rVert_{\mathcal{H}_2},
	\end{equation}
	where
	$
	{\eta}_i(s) : =(\mathbf{e}_i^\top \otimes C) \left[ M \otimes (sI_\ell - A) + L \otimes BC \right](F \otimes B),
	$ if \eqref{sysh} is considered, and 
	$
	{\eta}_i(s) : = \mathbf{e}_i^\top (sM + L )^{-1} F,
	$ if \eqref{sys} is considered.
\end{defn}
The transfer matrix $\eta_i(s)$ is the mapping from the external control signal $u$ to the output of the $i$-th subsystem, $y_i$, and thus $\eta_i(s)$ is interpreted as the behavior of the $i$-th vertex with respect to the external inputs. The concept of dissimilarity indicates how different two vertices are in terms of their behaviors. 
It is verified in \cite{ChengTAC2018MAS} that if the network system \eqref{sysh} is synchronized, $\mathcal{D}_{ij}$ in \eqref{eq:dissim} is well-defined, and a dissimilarity matrix $\mathcal{D} \in \mathbb{R}^{n \times n}$, with $[\mathcal{D}]_{ij} = \mathcal{D}_{ij}$ is symmetric, and with zero diagonal elements and nonnegative off-diagonal entries. However, it could be a formidable task to compute the dissimilarity between each pair of vertices in a large-scale network based on its definition. Next, we discuss efficient methods for computing dissimilarity $\mathcal{D}_{ij}$.

First, we consider the single-integrator network in \eqref{sys}, which is a semistable system. Following Section~\ref{sec:semistable}, the pseudo controllability Gramian of \eqref{sys} is computed as $\mathcal{P} = \mathcal{J} \tilde{\mathcal{P}} \mathcal{J}^\top$, where $\tilde{\mathcal{P}}$ is an arbitrary solution of
\begin{equation} \label{eq:J}
- M^{-1} L \tilde{\mathcal{P}} - \tilde{\mathcal{P}} L M^{-1}+ (I - \mathcal{J}) M^{-1} F F^{\top} M^{-1} (I - \mathcal{J}) = 0, \quad \mathcal{J}: = \frac{\mathds{1} \mathds{1}^\top M}{\mathds{1}^\top M \mathds{1}}.
\end{equation} 
We refer to
\cite{ChengCDC2016Gramian,ChengAuto2020Gramian} for more details. Note that Theorem~\ref{thm:H2semi} implies that the error system $S - \hat{S}$, with $S$, $\hat{S}$ defined in \eqref{eq:TFS}, is in the $\mathcal{H}_2$ space, and an   efficient method for computing  $\mathcal{D}$ is presented based on the pseudo controllability Gramian:
\begin{equation} \label{eq:Dij1}
\mathcal{D}_{ij} = \sqrt{(\mathbf{e}_i - \mathbf{e}_j) \mathcal{P} (\mathbf{e}_i - \mathbf{e}_j)^\top}.
\end{equation} 

Next, we consider the network system \eqref{sysh} which achieves synchronization. If the overall system \eqref{sysh} is semistable, we can still apply pseudo Gramians to compute dissimilarity. However, the subsystems in the network may be unstable. In this case, we present another computation approach \cite{ChengTAC2018MAS}. Denote 
\begin{equation}
\mathcal{S}: = 
\begin{bmatrix}
-I_{n-1} \\ \mathds{1}_{n-1}^\top 
\end{bmatrix}  
\in \mathbb{R}^{{n} \times {(n-1)}}, 
\quad
\mathcal{S}^\dagger = (\mathcal S^\top M^{-1}\mathcal S)^{-1}  \mathcal S^\top M^{-1},
\end{equation}
which satisfy $\mathcal{S} \mathds{1} = 0$, $\mathcal{S}^\dagger M  \mathds{1} = 0$, and $\mathcal{S}^\dagger \mathcal{S} = I_{n-1}$. Let 
\begin{equation*}
\mathcal{A}: = I_{n-1} \otimes A  -  \mathcal{S}^\dagger L M^{-1} {\mathcal{S}} \otimes BC, \quad \mathcal{B} = \mathcal{S}^\dagger F \otimes B,
\end{equation*}
where $\mathcal{A}$ is Hurwitz if and only if the system \eqref{sysh} achieves synchronization. 
\\

\begin{thm}
	\label{thm:clustercomput}
	Let the network system \eqref{sysh} achieve synchronization. Then, there exists a symmetric matrix
	$\bar{\mathcal{P}} \in \mathbb{R}^{(n-1) \ell \times (n-1)\ell}$, which is the unique solution of the Lyapunov equation
	$
	\bar{A} \bar{\mathcal{P}} + \bar{\mathcal{P}} \bar{A} + \bar{\mathcal{B}} \bar{\mathcal{B}}^\top = 0.
	$
	Moreover, 
	\begin{equation} \label{eq:Dij2}
	\mathcal{D}_{ij} = \sqrt{\tr(\Psi_{ij} \bar{\mathcal{P}} \Psi_{ij}^\top)},
	\end{equation}
	where $\Psi_{ij}: = (\mathbf{e}_i-\mathbf{e}_j)^\top M \mathcal{S} \otimes C$.
\end{thm}

The definition of pairwise dissimilarity in \eqref{eq:dissim} measures how close two subsystems behave, and aggregating vertices with similar behaviors potentially leads to a small approximation error. Having dissimilarity as a metric, clustering algorithms for static graphs in e.g., \cite{SurveyClustering,Schaeffer2007SurveyClustering} can be also adopted to solve the model reduction problem for dynamical networks. For instant, a \textit{hierarchical clustering} algorithm is applied in \cite{ChengTAC20172OROM} as follows.
\begin{algorithm} 
	\caption{Hierarchical Clustering Algorithm}
	\begin{algorithmic}[1]
		%        \Input $M$, $L$ and $B$, 
		%        model order $n$, desired order $r$
		%        \Output $\hat{M}$, $\hat{L}$, and $\hat{B}$ 
		\State Compute the dissimilarity matrix $\mathcal{D}$.
		\State Place each node into a singleton cluster, i.e., 
		$\mathcal{C}_i \leftarrow \{i\}, \ \forall \ 1 \leq i \leq n.$
		\State Find two clusters $\mathcal{C}_k$ and $\mathcal{C}_l$ such that
		\begin{equation} \label{eq:clusterDis}
		(k,l): =  \arg\min \left(    \dfrac{1}{|\mathcal{C}_k| \cdot |\mathcal{C}_l|}\sum_{i \in \mathcal{C}_k} \sum_{j \in \mathcal{C}_l}\mathcal{D}_{ij}\right). 
		\end{equation}
		\State Merge clusters $\mathcal{C}_k$ and $\mathcal{C}_l$ into a single cluster.
		\State Repeat the steps 3 and 4, until $r$ clusters are obtained. 
		
		\State Compute the characteristic matrix $\Pi \in \mathbb{R}^{n \times r}$ and return  
		$$\hat{M} \leftarrow \Pi^\top M \Pi , 
		\hat{L} \leftarrow \Pi^\top L \Pi ,  \hat{F} \leftarrow \Pi^\top F.$$
	\end{algorithmic}
	\label{alg3}
\end{algorithm}

An iterative approach for single-integrator networks can be found in \cite{ChengECC2016}, and an alternative K-means clustering method is presented in \cite{UmarCDC2019Clustering}, which takes into account the connectedness of vertices such that the vertices in each cluster form a connected graph.

In Algorithm \ref{alg3}, the proximity of two clusters $\mathcal{C}_\mu$  and $\mathcal{C}_\nu$ is evaluated by \eqref{eq:clusterDis}, which means the average dissimilarity of the vertices in the two clusters. 
%\begin{equation}  
%d(\mathcal{C}_\mu, \mathcal{C}_\nu) = \dfrac{1}{|\mathcal{C}_\mu| \cdot |\mathcal{C}_\nu|}\sum_{i \in \mathcal{C}_\mu} \sum_{j \in \mathcal{C}_\nu}\mathcal{D}_{ij}.
%\end{equation}
Other metrics of cluster proximity can be used as well. For instance, we can take the smallest dissimilarity of the vertices from two clusters, or the largest dissimilarity of the nodes from two clusters.
The proximity of two clusters allows us to link pairs of clusters with smaller proximity and place them into binary clusters. Then, the newly formed clusters can be grouped into larger ones according to the cluster proximity. In each loop, two clusters with the lowest proximity are merged together, and finally a binary hierarchical tree, called \textit{dendrogram}, that visualizes this process can be generated, see Fig.~\ref{example2_dendrogram} in the following example.

\begin{exm}
	\begin{figure}[!tp]\centering
		\includegraphics[scale=.46]{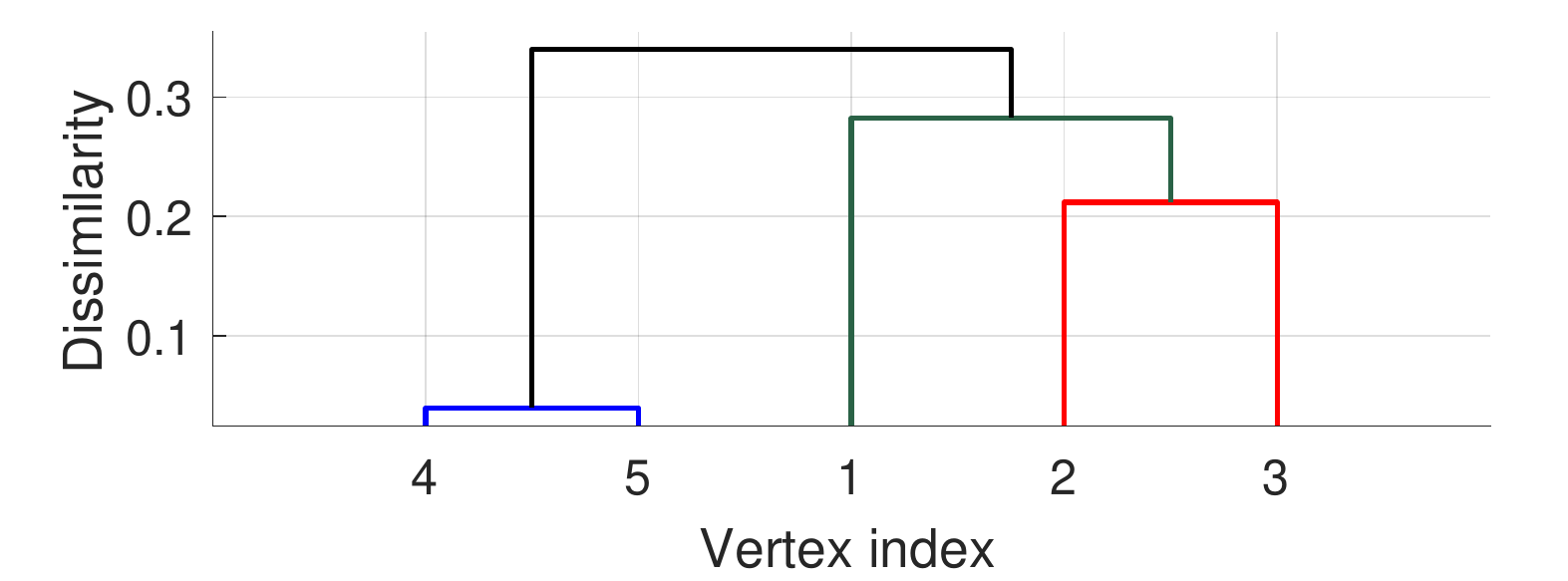}    
		\caption{Dendrogram illustrating the hierarchical clustering of the networked mass-damper system. The horizontal axis is
			labeled by vertex numberings, while the vertical axis represents
			the dissimilarity of clusters. The dissimilarity is measured in the $\mathcal{H}_2$-norm, and the level at which branches merge indicates the dissimilarity between two clusters.}
		\label{example2_dendrogram}
	\end{figure}    
	Consider the networked mass-damper system in Example \ref{exm1}.
	%    From Theorem \ref{thm:clustercomput}, we have 
	%    \begin{equation}
	%    \bar{\mathcal{A}} = \left[\begin{array}{rrrr}   
	%    -7  &  3 & -3  &  1 \\
	%    2  & -4 &  -2 & -1\\
	%    -1  &  1  & -9 &  1\\
	%    1   &      0  & -1  & -6
	%    \end{array}\right], \
	%    \bar{\mathcal{B}} = \left[\begin{array}{rr}   
	%    -1&   1\\
	%    0  & 1\\
	%    0   &1\\
	%    0   &0
	%    \end{array}\right]. 
	%    \end{equation}
	%    Solving the Lyapunov equation in \eqref{eq:LyapPro} yields
	%    \begin{equation}
	%    \bar{\mathcal{P}} = \left[\begin{array}{rrrr}   
	%    0.1716 &   0.1285  &  0.0647  &  0.0098\\
	%    0.1285  &  0.1476  &  0.0800 &   0.0066\\
	%    0.0647  &  0.0800  &  0.0573 &   0.0004\\
	%    0.0098 &   0.0066  &  0.0004 &   0.0016
	%    \end{array}\right]
	%    \end{equation}
	The dissimilarity matrix can be computed using either \eqref{eq:Dij1} or \eqref{eq:Dij2}, which yields
	\begin{equation*}  
	\mathcal{D}=\left[                
	\begin{array}{rrrrr}   
	0  &  0.2494  &  0.3154  &  0.3919  &  0.4142\\
	0.2494 &        0  &  0.2119  &  0.3688  &  0.3842\\
	0.3154 &   0.2119  &       0  &  0.2410  &  0.2394\\
	0.3919 &   0.3688  &  0.2410   &      0  &  0.0396\\
	0.4142 &   0.3842  &  0.2394  &  0.0396   &      0
	\end{array}
	\right].           
	\end{equation*}    
	The minimal value is $0.0396$, indicating that vertices 4 and 5 have the most similar behavior compared to the other pairs of vertices. Thus, vertices 4 and 5 are first aggregated, which leads to clusters: $\{ \{ 1\},\{ 2\},\{3\} \}, \{ 4,5\} \}$. In the hierarchical clustering, we check the proximities of the clusters by \eqref{eq:clusterDis} and then obtain a coarser clustering $\{ \{ 1\},\{ 2, 3\},\{ 4,5\} \}$. This process can be continued until we have generated a dendrogram as depicted in Fig. \ref{example2_dendrogram}.
\end{exm}

%\begin{rem} \label{rem:extadj}
%    Note that the formation of clusters in Algorithm \ref{alg3} does not focus on manipulating any individual edges. Even if two nodes are not adjacent, they can be placed into the same cluster if they have very similar behaviors.    
%    However, Algorithm \ref{alg3} can be easily adapted to only aggregate adjacent vertices. Define adjacent clusters as follows: Two clusters $\mathcal{C}_\mu$ and $\mathcal{C}_\nu$ are adjacent if there exist $i \in \mathcal{C}_\mu$ and $j \in \mathcal{C}_\nu$ such that vertices $i$ and $j$ are connected by an edge. To only aggregate adjacent vertices, we may add as an extra constraint to the step 3 in Algorithm \ref{alg3} that the clusters $\mathcal{C}_\mu$ and $\mathcal{C}_\nu$ are adjacent.
%\end{rem}

Algorithm~\ref{alg3} is based on pairwise dissimilarities of the vertices and minimizes within-cluster variances. The variance within a cluster can be characterized by the largest dissimilarity between all pairs of vertices within the cluster, which leads to an upper bound on the $\mathcal{H}_2$ approximation error \cite{ChengTAC2018MAS}.
% For simplicity, we denote 
%\begin{equation} \label{eq:Tr}
%\mathcal{S}_r: = 
%\frac{1}{r}\begin{bmatrix}
%\bm{1}_{r-1}\bm{1}_{r-1}^\top-rI_{r-1}
%\\ \bm{1}_{r-1}^\top
%\end{bmatrix}  
%\in \mathbb{R}^{r \times (r-1)},
%\end{equation}
%and 
%\begin{equation} \label{eq:MLbar}
%\bar{P} : = \Pi (\Pi^\top \Pi)^{-1} \mathcal{S}_r,
%\bar{M} := \bar{P}^\top \bar{P}, 
%\bar{L} := \bar{P}^\top L \bar{P}. 
%%    \bar{M} := \mathcal{S}_r \hat{M}^{-1} \mathcal{S}_r^\top, 
%%    \bar{L} := \mathcal{S}_r \hat{M}^{-1} \hat{L} \hat{M}^{-1} \mathcal{S}_r^\top,  
%\end{equation}
\vspace*{12pt}

\begin{thm} \label{thm:errbound}
	Consider the network system \eqref{sysh} with the output matrix $H = I_n$.
	%     Suppose $\bm{{\Sigma}}$ synchronizes (i.e., it satisfies Lemma \ref{lem:ConsCond1}. 
	Let  $\{\mathcal{C}_1, \mathcal{C}_2, ..., \mathcal{C}_r\}$ be the graph clustering of the network, and $G(s)$ and $\hat{G}(s)$ denote the transfer matrices defined in \eqref{eq:TF}. If $A$ in \eqref{sysagent} satisfies $A + A^\top < 0$, then it holds that
	\begin{equation} \label{eq:errbound}
	\lVert G(s) - \hat{G}(s) \rVert_{\mathcal{H}_2} <  \gamma \cdot \sum_{k=1}^{r} |\mathcal{C}_k| \cdot \max\limits_{i,j \in \mathcal{C}_k} \mathcal{D}_{ij},
	\end{equation} 
	where $\gamma  \in \mathbb{R}_{+}$ only depends on the original system \eqref{sysh} and satisfies 
	\begin{equation} \label{eq:LMIgamma}
	\begin{bmatrix}
	{I} \otimes (A^\top + A) - {L} \otimes (C^\top B^\top + BC) &   L \otimes BC & - I \otimes C^\top\\ 
	L \otimes C^\top B^\top & - \gamma  I & I \\
	- I\otimes C &  I & - \gamma  I 
	\end{bmatrix}   < 0.
	\end{equation}
	%    with $\mathcal{X}: = {I} \otimes (A^\top + A) - {L} \otimes (C^\top B^\top + BC)$. 
\end{thm}

If the considered network system is in the form of \eqref{sys}, we further obtain an error bound based on the pseudo controllability Gramian.
\\

\begin{pro} \cite{ChengAuto2020Gramian}
	Let $S$ and $\hat{S}$ in \eqref{eq:TFS} be the transfer matrices of \eqref{sys} and \eqref{sysr}, respectively. We have
	\begin{equation}
	\lVert S - \hat{S}\rVert_{\mathcal{H}_2} \leq 
	\gamma_s \sqrt{\tr(I-\Pi \Pi^\dagger) \mathcal{P}(I-\Pi \Pi^\dagger)^\top},
	\end{equation}
	where $\Pi^\dagger = (\Pi^\top M \Pi)^{-1} \Pi^\top M$, and $\mathcal{P}$ is the pseudo controllability Gramian of \eqref{sys}. The constant $\gamma_s \in  \mathbb{R}_{+}$ is a solution of
	\begin{equation}
	\begin{bmatrix}
	M L M^{-1} + M^{-1} L M & M^{-1} L & (I - \mathcal{J}^\top) H^\top \\
	L M^{-1} & -\gamma_s I & H^\top\\
	H (I - \mathcal{J}) & H & -\gamma_s I
	\end{bmatrix} \leq 0, 
	\end{equation}
	with $\mathcal{J}$ defined in \eqref{eq:J}.
\end{pro}
\vspace*{12pt}

The core step in dissimilarity-based clustering is to properly define the dissimilarity of dynamical vertices. For linear time-variant networks, nodal dissimilarity can be always defined as the transfer from the external inputs to the vertex states. This mechanism of dissimilarity-based clustering is applicable to different types of dynamical networks, see, e.g., \cite{ChengTAC20172OROM,ChengTAC2019Digraph,ChengTAC2019Digraph,ChengACOM2018Power} for more results on second-order networks, directed networks and controlled power networks. For nonlinear networks, 
DC gain, a function of input amplitude, can be considered \cite{kawano2019data}, in which model reduction is aggregating state variables having similar DC gains.  

\newpage 

\subsection{Edge Weighting Approach}

Generally, all the existing clustering-based reduction methods fall into the framework of Petrov-Galerkin Projections. In \cite{ChengECC2019weight,ChengTAC2020Weight}, an $\mathcal{H}_2$ optimal approach is presented, which does not aim to find a suitable graph clustering. Instead, this approach focuses on how to construct a ``good'' reduced-order model for a given clustering.
To formulate this problem, the topology of a reduced network can be obtained from the given clustering, while all the edge weights are free parameters to be determined via optimization algorithms.

Consider the original network system in \eqref{sys} with graph $\mathcal{G}$. Let $\{ \mathcal{C}_1,\mathcal{C}_2,\cdots,\mathcal{C}_r \}$ be a given graph clustering of $\mathcal{G}$. Then, a \textit{quotient graph} $\hat{\mathcal{G}}$ is a $r$-vertex directed graph obtained by aggregating all the vertices in each cluster as a single vertex, while retaining connections between clusters and ignoring the edges within each cluster. If there is an edge $(i,j) \in \mathcal G$ with vertices $i,j$ in the same cluster, then this edge will be ignored in $\hat{\mathcal{G}}$. However, if the edge $(i,j)$ satisfies $i \in \mathcal{C}_k$ and $j \in \mathcal{C}_l$, then there will be an edge $(k,l)$ in $\hat{\mathcal{G}}$. The incidence matrix $\hat{R}$ of the quotient graph $\hat{\mathcal{G}}$ can be obtained by removing all the zero columns of $\Pi^\top R$, where $R$ is the incidence matrix of $\mathcal G$, and $\Pi$ is the characteristic matrix of the clustering.
Denote 
\begin{equation}\label{omeg}
\hat{W} = \diag(\hat{w}), \ \text{with} \ \hat{w} = \begin{bmatrix}
\hat{w}_1 & \hat{w}_2 & \cdots & \hat{w}_\kappa 
\end{bmatrix}^\top, 
\end{equation} 
as the edge weight matrix of $\hat{\mathcal{G}}$, where $\hat{w}_k \in \mathbb{R}_{+}$, and $\kappa$ denotes the number of edges in $\hat{\mathcal{G}}$. Then, a parameterized model of a reduced-order network is obtained:
\begin{equation} \label{sysrp}
\left \{
\begin{array}{rcl}
\hat{M} \dot{z} &=& - \hat{R} \hat{W} \hat{R}^\top z +  \hat{F}u, \\
\hat{y} &=& \hat{H} z.
\end{array}
\right.
\end{equation} 
where $\hat{M} = \Pi^\top M \Pi$, $\hat{F} = \Pi^\top {F}$, and $\hat{H} = H \Pi$. The edge weight matrix $\hat{W}$ is the only unknown to be determined. Let 
\begin{equation} \label{eq:TFrp}
S_p : = \hat{H} (s \hat{M} + \hat{R} \hat{W} \hat{R}^\top)^{-1} \hat{F}.
\end{equation}
Then, an optimization problem can formulated to minimize the  approximation error $\lVert S - S_p \rVert_{\mathcal{H}_2}$ by tuning the edge weights. Here, an example is used to demonstrate the parametrized modeling of a reduced network system. \\

\begin{exm}
	Consider an undirected graph composed of 6 vertices in Fig.~\ref{fig:Example6Nodes}. An external force $u$ is acting on vertex 3, and the state of vertex 4 is measured as an the output $y$.  Given a clustering with 
	$\mathcal{C}_1 = \{1,2\},~ \mathcal{C}_2 = \{3\},~ \mathcal{C}_3 = \{4\},~ \mathcal{C}_4 = \{5,6\},$ the quotient graph is obtained in Fig.~\ref{fig:Example4Nodes} with the incidence matrix 
	\begin{equation*}
	\hat{\mathcal{B}} = \begin{bmatrix}
	1 & 1& 0 &0 \\
	-1 & 0& 1 &0\\
	0 & -1& 0 &1\\
	0 & 0& -1 &-1 \\
	\end{bmatrix}.
	\end{equation*}
	Let $\hat{W} = \diag(\hat{w}_1, \hat{w}_2, \hat{w}_3, \hat{w}_4)$ be the weights of the corresponding edges. The Laplacian matrix of the reduced network is constructed as
	\begin{equation*}
	\hat{R} \hat{W} \hat{R}^\top = \begin{bmatrix}
	\hat{w}_1 + \hat{w}_2 & - \hat{w}_1 & - \hat{w}_2 & 0\\
	-\hat{w}_1 & \hat{w}_1 + \hat{w}_4 & 0 & - \hat{w}_4\\
	- \hat{w}_2 & 0 & \hat{w}_2 + \hat{w}_3 & -\hat{w}_3\\
	0 & - \hat{w}_4 & - \hat{w}_3 & \hat{w}_3 + \hat{w}_4\\
	\end{bmatrix}, 
	\end{equation*}
	and moreover, we have $\hat{F} = \Pi^\top F =  \begin{bmatrix}
	0 & 1 & 0 & 0
	\end{bmatrix}^\top$, and $
	\hat{H} = H \Pi = \begin{bmatrix}
	0 & 0 & 1 & 0
	\end{bmatrix}$. If in the original network, $M = I_6$, in the reduced-order model \eqref{sysrp}, $\hat{M} = \Pi^\top M \Pi = \diag (2, 1, 1, 2)$.
	
	\begin{figure}[!tp]\centering
		\begin{minipage}[t]{0.5\linewidth}
			\centering
			\includegraphics[width=0.84\textwidth]{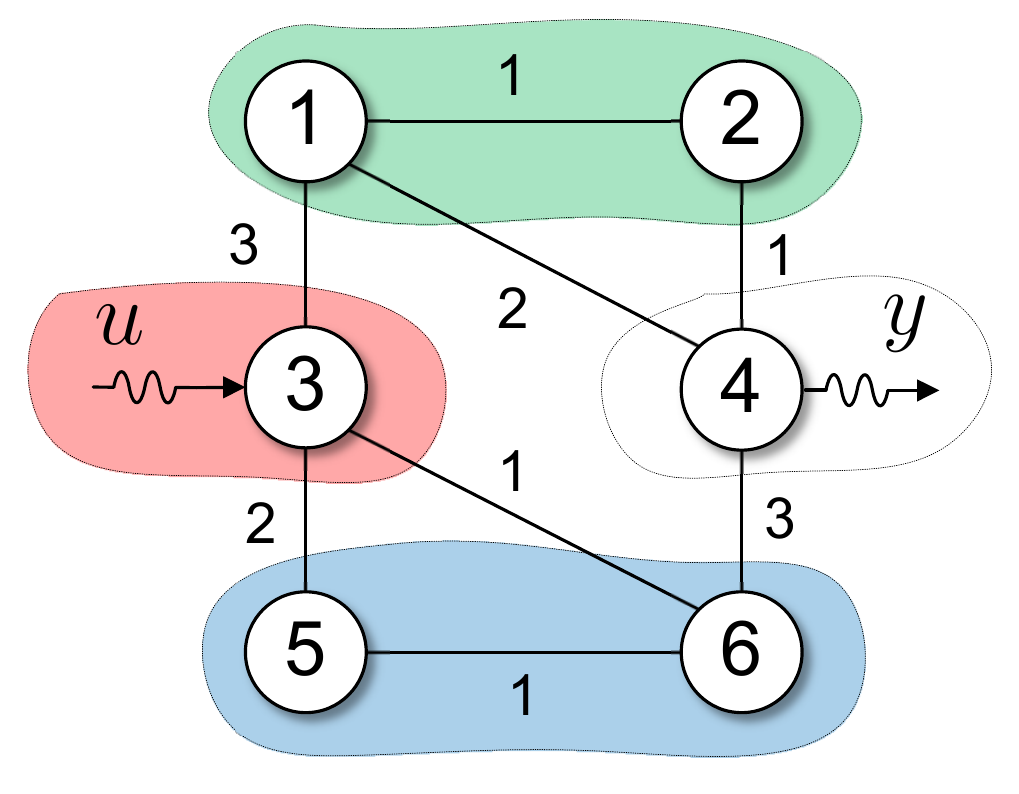}    
			\subcaption{}
			\label{fig:Example6Nodes}
		\end{minipage}%
		\begin{minipage}[t]{0.5\linewidth}
			\centering
			\includegraphics[width=0.75\textwidth]{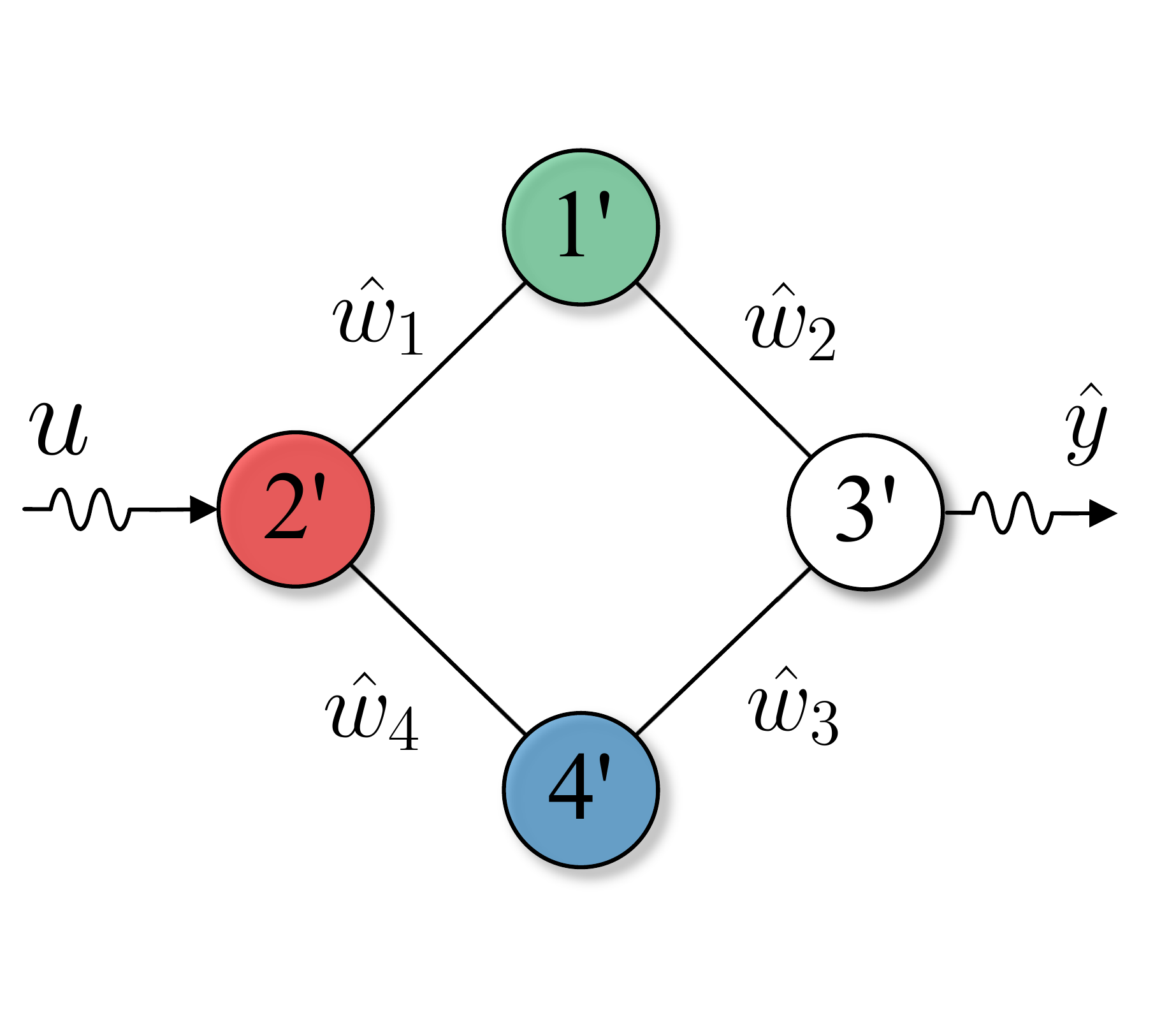}    
			\subcaption{}
			\label{fig:Example4Nodes}
		\end{minipage}%
		\caption{(a) An undirected network consisting of 6 vertices, in which vertex 3 is the leader and vertex 4 is measured. Four clusters are indicated by different colors. (b) A quotient graph obtained by clustering.}
	\end{figure}
\end{exm}

An optimization technique based on the \textit{convex-concave decomposition} can be applied to search for a set of optimal weights iteratively. Before proceeding, a necessary and sufficient condition for characterizing $\lVert G_e(s)\rVert_{\mathcal{H}_2}$ is shown.

\begin{thm}\label{theorem-1}
	Given the network system \eqref{sys}. A reduced-order model in \eqref{sysrp} satisfies $\|S - S_p\|_{\mathcal{H}_2}^2<\gamma_p$ if and only if there exist matrices  $\hat{Q}=\hat{Q}^\top>0$,
	$\hat{Z}=\hat{Z}^\top>0$ ,
	%  with dimension $R\in\mathbb{R}^{q\times q}$,
	and $\hat{\delta} \in \mathbb{R}_{+}$ such that $\tr(\hat{Z})<\gamma_p$, and 
	\begin{eqnarray}
	\label{inquali-1}
	\begin{bmatrix}
	\hat{Q}\bar{A}+\bar{A}^\top\hat{Q}& \hat{Q}B_e& \hat{Q}E\\
	B^\top_{e}\hat{Q} & -\hat{\delta}I& 0\\
	E^\top\hat{Q}& 0& 0
	\end{bmatrix}+\begin{bmatrix}
	-\bar{A}^\top_{r}\bar{A}_{r}& 0& \bar{A}^\top_{r}\\
	0& 0& 0\\
	\bar{A}_{r}& 0 & -I
	\end{bmatrix}<0,
	\\
	\begin{bmatrix}\label{inquali-2}
	\hat{Q} & \hat{\delta}C^\top_e\\
	\hat{\delta}C_e & \hat{Z} \\
	\end{bmatrix}>0, 
	\end{eqnarray}
	%%%%%%%%%%%%%%%%%%%%%%%%%%%%%
	where 
	\begin{align*} 
	\bar{A} & =\begin{bmatrix}
	-\mathcal{S}_n^+ L  M^{-1} \mathcal S_n& 0\\
	0& 0
	\end{bmatrix}, \
	\bar{A}_r =\begin{bmatrix}
	0& -\mathcal S_r^+ \hat{R} \hat{W} \hat{R}^\top\hat{M}^{-1} \mathcal S_r\\
	0& 0
	\end{bmatrix}, \
	B_e  =  \begin{bmatrix} \mathcal S_n^+ F  \\ \mathcal S_r^+ \hat{F} \end{bmatrix},
	\\
	C_e & =\begin{bmatrix} H M^{-1}\mathcal S_n &  -\hat{H} \hat{M}^{-1}\mathcal S_r \end{bmatrix},
	\
	E  =\begin{bmatrix}
	0& 0\\
	I& 0
	\end{bmatrix},
	\mathcal{S}_n = 
	\begin{bmatrix}
	-I_{n-1} \\ \mathds{1}_{n-1}^\top 
	\end{bmatrix} , \ 
	\mathcal{S}_r = 
	\begin{bmatrix}
	-I_{r-1} \\ \mathds{1}_{r-1}^\top 
	\end{bmatrix}. 
	\end{align*}
\end{thm}

Based on Theorem \ref{theorem-1}, the edge weighting problem is formulated as an minimization problem: 
\begin{align}\label{pro:min-1}
\min_{\hat{Q}>0, ~\hat{W}} \tr(Z), \quad \quad 
\text{s.t.\quad  \eqref{inquali-1} and  \eqref{inquali-2} hold}, 
\end{align}
where $\hat{Z}=\hat{\delta}Z$. Consider the matrix-valued mapping:
\begin{equation}\label{mapping}
\Phi(\hat{Q}, \hat{\delta}, \hat{W})=\psi(\hat{Q}, \hat{\delta})+\varphi(\hat{W}), 
\end{equation}
where 
\begin{align*}
\psi(\hat{Q}, \hat{\delta}) =\begin{bmatrix}
\hat{Q}\bar{A}+\bar{A}^\top\hat{Q}& \hat{Q}B_e& \hat{Q}E \\
B^\top_{e}\hat{Q} & -\hat{\delta}I& 0\\
E^\top\hat{Q}& 0& 0
\end{bmatrix}, \
\varphi(\hat{W}) =\begin{bmatrix} 
-\bar{A}^\top_{r}\bar{A}_{r}& 0& \bar{A}^\top_{r}\\
0& 0& 0\\
\bar{A}_{r}& 0 & -I
\end{bmatrix}. 
\end{align*}
Then, the pair $(\psi, -\varphi)$ is a psd-convex-concave decomposition of $\Phi$ \cite{ChengTAC2020Weight}. The bilinear matrix inequality \eqref{inquali-1} with the nonlinearity term $\bar{A}^\top_{r}\bar{A}_{r}$ can be handled using such a decomposition, which  
can linearize the optimization problem \eqref{pro:min-1} at a stationary point $\hat{W}$ \cite{dinh2011combining}. Rewrite  $\varphi(\hat{W})$ in \eqref{mapping} as $
\phi(\hat{w}) = \varphi(\hat{W}),
$ with $\hat{w} \in \mathbb{R}_{+}^\kappa$ defined in \eqref{omeg}.
%The derivative of the matrix-valued mapping $\phi(\mu)$ at $\mu$ is a linear
%mapping $D \phi:  \mathbb{R}_+^{\bar{m}} \rightarrow \mathcal{S}^{\ell}$, with $\ell = n+2r+p+pq-3$, which is defined as
%\begin{equation}\label{eq:deriva}
%D \phi(\mu) [h] =\sum_{i=1}^{\bar{m}} h_{i}\frac{\partial\phi}{\partial \mu_{i}}(\mu), \ \forall~h \in \mathbb{R}^{\bar{m}}.
%\end{equation}
Given a point $\hat{w}^{(k)}$, the linearized formulation of the problem \eqref{pro:min-1} at $\hat{w}^{(k)}$ is formulated as a \textit{convex} problem:
\begin{align}\label{pro:lin-min1}
&\min_{\hat{Q}>0, \hat{w} \in \mathbb{R}_{+}^{\kappa} } f(\hat{w})=\tr(Z)
\\
&\quad \text{s.t.} \quad     \begin{bmatrix}
\hat{Q} & \hat{\delta}C^\top_e\\
\hat{\delta}C_e & \hat{Z} \\
\end{bmatrix}>0,\ \hat{\delta} \in \mathbb{R}_{+}, \ \hat{Z} = \hat{\delta} Z >0   \nonumber
\\
& \hspace*{35pt} \psi(\hat{Q}, \hat{\delta})+\varphi(\hat{W}^{(k)})+D\phi(\hat{w}^{(k)})[\hat{w}-\hat{w}^{(k)}]<0, \nonumber
\end{align}
where the derivative of $\phi(\hat{w}^{(k)})$ is defined as
$$
D\phi(\hat{w}^{(k)})[\hat{w}-\hat{w}^{(k)}]:=\sum_{i=1}^{\kappa}(\hat{w}_{i}-\hat{w}^{(k)}_{i})\frac{\partial\phi}{\partial \hat{w}_{i}^{(k)}}(\hat{w}^{(k)}).
$$
Then, an algorithmic approach is presented in Algorithm~\ref{alg} for
solving the minimization problem in \eqref{pro:min-1} in an iterative fashion.

\begin{algorithm} 
	\caption{Iterative Edge Weighting}
	\begin{algorithmic}[1]
		\State Choose an initial vector $\hat{w}^{(0)} \in \mathbb{R}_{+}^{\kappa}$.

		\State Set iteration step: $k \leftarrow 0$.
		\Repeat
		\State Solve \eqref{pro:lin-min1} to obtain the optimal solution $\hat{w}^{*}$.
		
		\State   $k \leftarrow k+1$, and $\hat{w}^{(k)} \leftarrow \hat{w}^{*}$.
		
		\Until{$|f(\mu^{(k+1)})-f(\mu^{(k)})| \leq \varepsilon$, with $\varepsilon$ a prefixed error tolerance.
			
			\State Return 
			$\hat{W}^{*} \leftarrow \diag(\hat{w}^*)$.}
	\end{algorithmic}
	\label{alg}
\end{algorithm}

If $\hat{w}$ is initialized as the outcome of clustering-based projection methods, the approximation accuracy obtained by edge weighting approach will be better than the ones obtained by clustering-based projection after iteration.  Furthermore, to solve the optimization problem in \eqref{pro:min-1}, we can also use a cross iteration algorithm presented in \cite{ChengECC2019weight}.

\subsection{Other Clustering-Based Methods}
\label{XD.sec:othercluster}

In this section, several other model reduction schemes based on graph clustering are reviewed. 
The method in \cite{ChengECC2020} formulates the clustering-based model reduction as a nonconvex optimization problem with mixed-binary
variables $\Pi$ and the objective function to minimize the $\mathcal{H}_2$-norm of the approximation error. The error system between \eqref{sys} and \eqref{sysr} is defined, of which the controllability and the observability Gramians are used to derive an explicit expression for the gradient of
the objective function. Then a projected gradient algorithm can be employed to solve this optimization problem with mathematical guarantees on its convergence. Related to the work in \cite{Monshizadeh2014}, a combination of the Krylov subspace method with graph clustering is proposed in \cite{Petar2015CDC}, where a reduced basis is firstly found by the Iterative Rational Krylov Algorithm (IRKA), and then a graph partition is obtained by the QR decomposition with column pivoting on the projection matrix. 
%However, the meaning of the obtained clustering is not clear, and there is no error bound for the approximation. 
An alternative graph-based model reduction method is proposed in \cite{leiter2017graph}, which finds a graph clustering based on the edge agreement protocol of a network (see the definition in \cite{zelazo2013performance}) and provides a greedy contraction algorithm as a suboptimal solution of graph clustering.  The clustering and aggregation approach in \cite{deng2011optimal,deng2014structure} is 
based on spectral analysis of Markov chains
The Kullback–Leibler divergence rate is employed as a metric to measure the difference between the original network model and its approximation.

Clustering-based model reduction approaches are also found in the applications of other types of networks, i.e., network systems that do not reach consensus. Instead, other network properties are emphasized. For instance, \cite{martin2018large} proposes a reduction method for scale-free networks, which are networks whose degree distribution follows a power law. They are
roughly characterized by the presence of few vertices
with a large degree (number of connections) and a large number
of verities with small degree. The method in \cite{martin2018large} preserves the eigenvector centrality of the adjacency matrix of the original network such that the obtained reduced network remains scale-free.

Positive networks are considered in \cite{Ishizaki2014,ishizaki2015clustereddirected}. A single-input bidirectional positive network is given in \cite{Ishizaki2014} as
\begin{equation}\label{positivenet}
\dot{x} =  A x + b u, \ x \in \mathbb{R}^n, u \in \mathbb{R},
\end{equation}
where $b \in \mathbb{R}^p$, and  $A := - D -L$ with $D \geq 0$ a diagonal matrix (i.e., at least one diagonal entry of $D$ should be positive and the rest of the diagonal entries are zero) and $L \geq 0$ a Laplacian matrix representing an undirected connected graph. It is verified that $A$ is negative definite, and thus the system \eqref{positivenet} is asymptotically stable. The structure of $A$ can be interpreted as a network containing self-loops.  In \cite{Ishizaki2014}, a set of clusters is constructed based on the notion of cluster reducibility, which characterizes the uncontrollability of local state variables. By aggregating
the reducible clusters,
a reduced order model is obtained that preserves the stability and positivity. The work in \cite{ishizaki2015clustereddirected} extends this method to the directed case, where $A$ in \eqref{positivenet} is now assumed to be irreducible, Metzler and semistable. In this case, the Frobenius eigenvector of $A$ is used for constructing the projections such that both
semistability and positivity are preserved in the resulting reduced-order network model. In both \cite{Ishizaki2014} and \cite{ishizaki2015clustereddirected}, an upper bound on the approximation error is established
using the cluster reducibility, and then a clustering scheme is proposed to select suitable clusters, aiming at minimizing the \textit{posteriori}  bound on the reduction error.

%%%%%%%%%%%%%%%%%%%%%%%%%%%%%%%%%%%%%%%%%%%%%%%%%%%%%%%%%%%%%%%%%%%%%%%%%%%%%%%%

\section{Balanced Truncation of Network Systems}
\label{sec:BTofNetwork}
%%%%%%%%%%%%%%%%%%%%%%%%%%%%%%%%%%%%%%%%%%%%%%%%%%%%%%%%%%%%%%

Reducing the dimension of each subsystem also results in a simplification of overall networks. To reduce the dynamics of vertices, balanced truncation based on generalized Gramian matrices is commonly used, see, e.g.,   \cite{Monshizadeh2013stability,ChengEJC2018Lure,ChengAuto2019}, in which preserving the synchronization property of the overall network is of particular interest. In this section, we review some recent results in the synchronization-preserving model reduction of large-scale network system using the classic generalized balanced truncation. For simplicity, we assume $M = I_n$ in \eqref{sysh} throughout this section.

\subsection{Model Reduction of Subsystems in Networks}

Starting from a synchronized network system in \eqref{sysh}, the aim of this subsection is to derive a network model with reduced-order subsystems such that synchronization is preserved in the reduced-order network in \eqref{sysrh}. 

If each subsystem in \eqref{sysagent} is asymptotically stable, we might apply standard balanced truncation to reduce the dimension of the subsystem regardless of their interconnection structure. However, this reduction is possible to destroy the property of the overall network system \eqref{sysh}, e.g., stability and synchronization. To achieve synchronization preservation, \cite{Monshizadeh2013stability} adopts a sufficient small gain type of condition to guarantee synchronization of \eqref{sysh}.

\begin{lem} \label{lem:Riccati}
	Denote by $0 = \lambda_1 < \lambda_2 \leq \cdots \leq \lambda_{n}$ the eigenvalues of the Laplacian matrix $L$. 
	The network system \eqref{sysh} achieves synchronization if there exists a nonzero eigenvalue $\lambda \in \{\lambda_2, ..., \lambda_{n} \}$ such that $A - \lambda BC$ is Hurwitz, 
	and there exists a positive definite matrix K satisfying the Riccati inequality 
	\begin{equation} \label{eq:Riccati}
		(A - \lambda BC)^\top K + K(A - \lambda BC) + C^\top C + \left( \frac{\delta}{\gamma}\right)^2 K BB^\top K < 0,
	\end{equation}
	where $\delta: = \max \{\lambda - \lambda_2, \lambda_n - \lambda \}$.
\end{lem}

It is worth noting that \eqref{eq:Riccati} is equivalent to the small gain condition $$\| C (sI_\ell - A + \lambda BC)^{-1} B \|_{\mathcal{H}_\infty} < \frac{\delta}{\gamma}.$$ 
Let $K_m$ and $K_M$ be the minimal and maximal real symmetric solutions of \eqref{eq:Riccati}, then $K_M^{-1}$ and $K_m$ can be regarded as a pair of generalized Gramians of the system $(A + \lambda BC, \frac{\delta}{\gamma} B, C)$. Apply the generalized balanced truncation introduced in Section~\ref{sec:PreBT}, a reduced-order model $(\hat{A} + \lambda \hat{B}\hat{C}, \frac{\delta}{\gamma} \hat{B}, \hat{C})$ with $\hat{A} \in \mathbb{R}^{k \times k}$ is obtained such that the small gain condition 
$\| \hat{C} (sI_\ell - \hat{A} + \lambda \hat{B}\hat{C})^{-1} \hat{B} \|_{\mathcal{H}_\infty} < \frac{\delta}{\gamma}$ is retained. Therefore, the following reduced-order network model preserves the synchronization property. 
\begin{equation} \label{syshsubred}
\left\{
\begin{array}{lr}
\dot{\xi} = (I_n \otimes \hat{A} - L \otimes \hat{B}\hat{C}) \xi  + (F \otimes \hat{B}) u, \\
\eta = 
(H \otimes \hat{C}) \xi.\\	
\end{array}
\right.
\end{equation}

\begin{thm}
	Consider a network system \eqref{sysh} that satisfies the synchronization condition in Lemma~\ref{lem:Riccati}. Then, the reduced-order network model in \eqref{syshsubred}
	obtained by generalized balanced truncation using $K_M^{-1}$ and $K_m$ achieves synchronization.
\end{thm}

Moreover, similar to \eqref{eq:outputAEP}, we assume a particular output  
$
y = (W^{\frac{1}{2}} R^\top \otimes C) x,
$
then error system between \eqref{sysh} and \eqref{syshsubred} is stable. Denote $\tilde{G}(s)$ as the transfer matrix of system \eqref{syshsubred}, and the model reduction error is upper bounded as 
\begin{equation}
	\| G(s) - \tilde{G} (s)\|_{\mathcal{H}_\infty} \leq  \frac{2 \gamma \sqrt{\lambda_{n}}}{\delta (1 - \gamma^2)}  \sum_{i = k+1}^{\ell} \sigma_i
\end{equation}
where $\sigma_i$ are the GHSVs computed using $K_M^{-1}$ and $K_m$  \cite{Monshizadeh2013stability}.

Inspired by the work \cite{Monshizadeh2013stability} for linear networks,  \cite{ChengECC2018Lure,ChengEJC2018Lure} considers dynamical networks of diffusively interconnected nonlinear Lur'e subsystems. The robust synchronization of the Lur'e network can be characterized by a linear matrix inequality (LMI). Different from \cite{Monshizadeh2013stability,ChengECC2018Lure}, where the minimum and maximum solutions of the LMI are used as generalized Gramians, \cite{ChengEJC2018Lure} suggests to only use the solution of the LMI with the minimal trace as a generalized controllability Gramian, while the observability counterpart is taken by the standard observability Gramian as the solution of Lyapunov equation, which is less conservative than the LMI. Using such a pair, generalized balanced truncation is performed on the linear component of each Lur'e subsystem, and the resulting reduced order network system is still guaranteed to have the robust synchronization property. 

\subsection{Simultaneously Reduction of Network Structure and Subsystems}

In the line of works \cite{Ishizaki2014,ishizaki2015clustereddirected}, a reduction method for network systems
composed of higher-dimensional dissipative
subsystems is presented in \cite{Ishizaki2016dissipative}, where the subsystems are reduced via block-diagonally orthogonal projection, while the network structure is simplified using clustering. In \cite{ChengIFAC2017}, the balancing method, for the first time, is applied for reducing the interconnection structure of networks with diffusively coupled vertices, and more extensions are found in	\cite{LanlinCDC2019eigenvalue,LanlinCDC2019second-order} based on eigenvalue assignment and moment matching.
In \cite{ChengAuto2019}, the idea in \cite{ChengIFAC2017} is further developed and applied to general networks of the form \eqref{sysh}. The proposed approach can reduce the
complexity of network structures and individual agent
dynamics simultaneously via a unified framework.  

%\subsection{Separation of Network System}
%\label{sec:Separation}
Consider the network system $\bm{\Sigma}$ in \eqref{sysh}, where each subsystem $\bm{\Sigma}_i$ as in \eqref{sysagent} is passive, namely, there exists a positive definite $K$ such that \eqref{eq:passive} holds. Note that $L$ is singular and $A$ in \eqref{sysagent} is not necessarily Hurwitz, implying that the overall system $\bm{\Sigma}$ may be not asymptotically stable, and thus a direct application of balanced truncation to $\bm{\Sigma}$ is not feasible. The method in \cite{ChengAuto2019} starts off with a decomposition of $\bm{\Sigma}$ using a \textit{spectral decomposition} of the graph Laplacian:  
\begin{equation} \label{eq:LEigDec}
L = T {\it \Lambda} T^\top = \begin{bmatrix}
T_1 & T_2
\end{bmatrix} 
\begin{bmatrix}
\ibLambda & \\  & 0
\end{bmatrix}
\begin{bmatrix}
T_1^\top \\ T_2^\top
\end{bmatrix},
\end{equation}
where $T_2 = \mathds{1}_n/\sqrt{n}$ and
$
\ibLambda : = \text{diag}(\lambda_2, \lambda_3, ..., \lambda_{n}),
$
with $\lambda_i$ the nonzero eigenvalues of $L$. Then, the system $\bm{\Sigma}$ can be split into two   components, namely, an \textit{average module}
\begin{equation} \label{sysa}
{\bf{\Sigma_a}}: 
\left\{
\begin{array}{lr}
\dot{z}_a = Az_a + \dfrac{1}{\sqrt{n}}  (\mathds{1}_n^\top F\otimes B) u, \\
y_a = \dfrac{1}{\sqrt{n}}  ( H \mathds{1}_n \otimes C) z_a, \\	
\end{array}
\right.
\end{equation}
with $z_a \in \mathbb{R}^{\ell}$,
and an asymptotically stable system
\begin{equation} \label{syst}
{\bf{\Sigma_s}}: 
\left\{
\begin{array}{lr}
\dot{z}_s = (I_{N-1} \otimes A - \ibLambda \otimes BC) z_s + ( \bar{F} \otimes B) u, \\
y_s = (\bar{H} \otimes C) z_s.\\	
\end{array}
\right.
\end{equation}
where $z_s \in \mathbb{R}^{(n-1)\times \ell}$, $\bar{F}=T_1^{T} F$, {and} $\bar{H} = H T_1$. The stability is guaranteed as the system $\bm{\Sigma}$ achieves synchronization Theorem~\ref{thm:syshsyn}.

\begin{figure}[!tp]\centering
	\centering
	\includegraphics[width=0.75\textwidth]{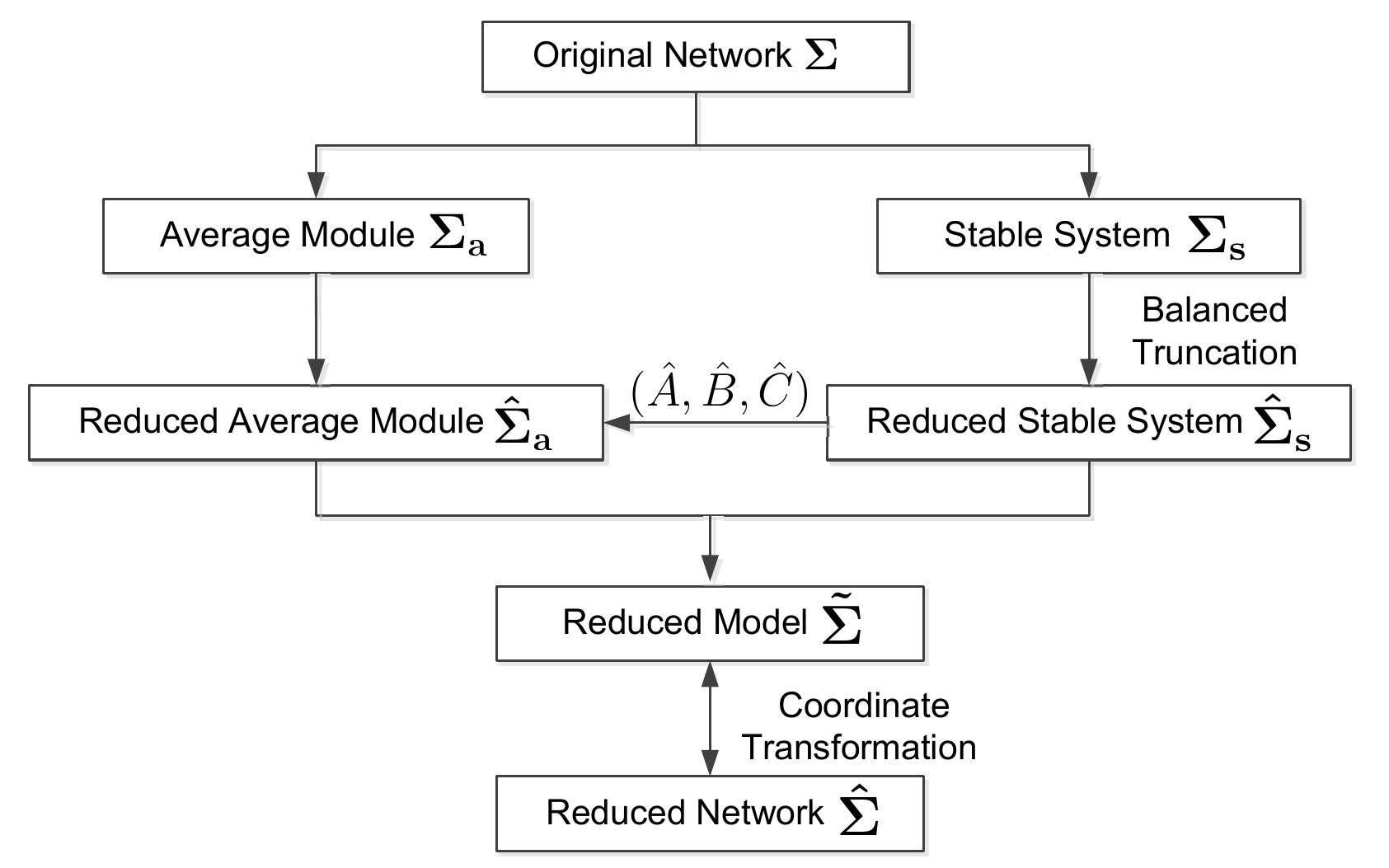}	
	\caption{The model reduction scheme for networked passive systems, where the simplification of network structure and the reduction of subsystems are performed simultaneously.}	
	\label{procedure}
\end{figure}

The model reduction procedure is as follows. First, we can apply
balanced truncation to $\bf{{\Sigma}_s}$ to generate a lower-order approximation $\bf{\hat{\Sigma}_s}$. It meanwhile gives a reduced subsystem  $(\hat{A}, \hat{B}, \hat{C})$ resulting in a reduced order average module $\bf{\hat{\Sigma}_a}$. Combining  $\bf{\hat{\Sigma}_s}$ with $\bf{\hat{\Sigma}_a}$ then formulates a reduced order model $\bm{\tilde{\Sigma}}$ whose input-output behavior approximates that of the original system $\bm{\Sigma}$. However, at this stage, the network structure is not necessarily preserved by $\bm{\tilde{\Sigma}}$. Then, it is desired to use a coordinate transformation to convert $\bm{\tilde{\Sigma}}$ to $\hat{\bm{\Sigma}}$, which restores the Laplacian structure. The whole procedure is summarized in Fig. \ref{procedure}. There are two key problems here: 
\begin{enumerate}
	\item How to retain the subsystem structure in $\bf{\hat{\Sigma}_s}$ such that subsystem dynamics do not mix with the topological information?
	
	\item How to recover a network interpretation in the reduced-order model $\bm{\tilde{\Sigma}}$ via a coordinate transformation?
\end{enumerate} 

To resolve the first problem, we resort to the balanced truncation approach based on generalized Gramians.
Suppose $\ibLambda$ in \eqref{eq:LEigDec} has $s$ distinct diagonal entries ordered as: $\bar{\lambda}_1 > \bar{\lambda}_2 > \cdots > \bar{\lambda}_s$. We rewrite $\ibLambda$ as
$ 
\ibLambda = \mathrm{blkdiag} (\bar{\lambda}_1 I_{m_1}, \bar{\lambda}_2 I_{m_2}, ... , \bar{\lambda}_s I_{m_s}),
$
where $m_i$ is the multiplicity of $\bar{\lambda}_i$, and $\sum_{i=1}^{s} m_i = n-1$. 
%In order to guarantee that the reduced order model will satisfy the desired properties in Problem \ref{prob:appx}, 
Then,    
%define the  generalized controllability and observability Gramians of \eqref{sys:accompany} as the solutions $X$ and $Y$ to 
the following Lyapunov equation and inequality have solutions $X$ and $Y$.
\begin{subequations} \label{eq:Lyap}
	\begin{align}
	-\ibLambda {X}  - {X}  \ibLambda + \bar{F}\bar{F}^\top &= 0, \label{eq:LyapX}
	\\ 
	-\ibLambda {Y}  - {Y} \ibLambda + \bar{H}^\top\bar{H} &\leq 0, \label{eq:LyapY} 
	\end{align}
\end{subequations}
where, ${X} = {X}^\top > 0$ and  
$  
{Y} := \blkdiag (Y_1, Y_2, ..., Y_s),
$
with $Y_i=Y_i^\top> 0$ and $Y_i \in \mathbb{R}^{m_i \times m_i}$, for $i = 1,2, ..., s$. 
The generalized controllability and observability Gramians of the stable system $\bf{{\Sigma}_s}$ are characterized by the following theorem.
\begin{thm} \label{thm:GenGram}
	Let $X>0$ be the unique solution of \eqref{eq:LyapX}, and let $Y > 0$ be a solution of \eqref{eq:LyapY}. Let $K_{m}>0$ and $K_{M}>0$ be the minimum and maximum solutions of \eqref{eq:passive}, respectively. 
	Then, the matrices
	\begin{equation} \label{eq:GenGrams}
	\mathcal{X} := {X} \otimes K_{M}^{-1} \ \text{and} \ \mathcal{Y}: =  {Y} \otimes K_m
	\end{equation} 
	are a pair of generalized Gramians of the asymptotically stable system $\bf{\Sigma_s}$.
	%	 i.e., satisfying the inequalities in \eqref{eq:ConGramInq} and \eqref{eq:ObsGramInq}, respectively. 	
	Moreover, there exist two nonsingular matrices $T_\mathcal{G}$ and $T_\mathcal{D}$ such that $\mathcal{T} = T_\mathcal{G} \otimes T_\mathcal{D}$ satisfies
	\begin{equation} \label{eq:balancing}
	\mathcal{T} \mathcal{X} \mathcal{T}^\top = \mathcal{T}^{-T}\mathcal{Y}\mathcal{T}^{-1} = \Sigma_\mathcal{G} \otimes \Sigma_\mathcal{D}.
	\end{equation}
	Here, $\Sigma_\mathcal{G} := \mathrm{diag} \{\sigma_1, \sigma_2,..., \sigma_{n-1}\},$ and $\Sigma_\mathcal{D} := \mathrm{diag} \{\tau_1, \tau_2,..., \tau_{\ell}\}$, where $\sigma_1 \geq \sigma_2 \geq \cdots \geq \sigma_{n-1},$ and $\tau_1 \geq \tau_2 \geq \cdots \geq \tau_{\ell}
	$ are
	equal to the square roots of the eigenvalues of $XY$ and $K_{M}^{-1}  K_{m}$, respectively. 
\end{thm}

The block-diagonal structure of $Y$ will be crucial to guarantee  that the reduced order model, obtained by preforming balanced truncation on the basis of $X$ and $Y$, can be interpreted as a network system again, as will be shown in Theorem~\ref{thm:LapReal} below.

\begin{rem} \label{rem:XY}
	By the duality between controllability and observability, we can also use $-\ibLambda {X}  - {X}  \ibLambda + \bar{F}\bar{F}^\top \leq  0$, and $-\ibLambda {Y}  - {Y} \ibLambda + \bar{H}^\top\bar{H} = 0$ to characterize the pair $X$ and $Y$ for the balanced truncation, where now $X$ is constrained to have a block-diagonal structure.
\end{rem}

Selecting the pair of Gramians in \eqref{eq:GenGrams} with the Kronecker product structure is meaningful, since they can be simultaneously diagonalized, (i.e., balanced) using transformations of the form $\mathcal{T} = T_\mathcal{G} \otimes T_\mathcal{D}$. Note that $T_\mathcal{G}$ and $T_\mathcal{D}$   are independently generated from \eqref{eq:Lyap} and \eqref{eq:passive}. More precisely, $T_\mathcal{G}$ only balances the network structure, or the triplet $(\ibLambda, \bar{F}, \bar{H})$, while $T_\mathcal{D}$ only balances the  agent dynamics, i.e., the triplet $(A, B, C)$. Thus, the Laplacian dynamics and each subsystem \eqref{sysagent} can be reduced independently, allowing the resulting reduced order model to preserve a network interpretation as well as the passivity of subsystems. 

Denote
$
	(\ihLambda_{1}, \hat{F}_{1}, \hat{H}_{1})\ \text{and}\  (\hat{A}, \hat{B}, \hat{C}): = {\bm{\hat{\Sigma}}_i}
$
	as the reduced order models of $(\ibLambda, \bar{F}, \bar{H})$ and $(A, B, C)$, respectively, where 
	$\ihLambda_{1} \in \mathbb{R}^{(r-1) \times (r-1)}$,
	$\hat{F}_{1} \in \mathbb{R}^{(r-1) \times  p}$, 
	$\hat{H}_{1} \in \mathbb{R}^{q \times (r-1)}$, $\hat{A} \in \mathbb{R}^{k \times k}$,
	$\hat{B} \in \mathbb{R}^{k \times  m}$, and
	$\hat{C} \in \mathbb{R}^{m \times k}$. Here, $k < \ell$ and $r < n$. 
Consequently, the reduced order models of the average module \eqref{sysa} and the stable system \eqref{syst} are constructed.
\begin{subequations}  
	\begin{align} 
	\label{sysaRed}
	\bf{\hat{\Sigma}_a}: &
	\left\{
	\begin{array}{lr}
	\dot{\hat{z}}_a = \hat{A} \hat{z}_a + \dfrac{1}{\sqrt{n}}  ( \mathds{1}_n^\top F\otimes \hat{B}) u, \\
	\hat{y}_a = \dfrac{1}{\sqrt{n}}  ( H  \mathds{1}_n \otimes \hat{C}) \hat{z}_a. \\	
	\end{array}
	\right.
	\\
	\label{systRed}
	\bf{\hat{\Sigma}_s}: &
	\left\{
	\begin{array}{lr}
	\dot{\hat{z}}_s = (I_{r-1} \otimes \hat{A} - \ihLambda_{1} \otimes \hat{B}\hat{C}) \hat{z}_s  + (\hat{F}_{1} \otimes \hat{B}) u, \\
	\hat{y}_s = 
	(\hat{H}_{1} \otimes \hat{C}) \hat{z}_s.\\	
	\end{array}
	\right.
	\end{align}
\end{subequations}  
Combining the reduced order models  $\bf{\hat{\Sigma}_a}$ and $\bf{\hat{\Sigma}_s}$,  a lower-dimensional approximation of the overall system $\bm{\Sigma}$ is formulated as 
\begin{equation} \label{sysRedMod}
\bm{\tilde{\Sigma}}: 
\left\{
\begin{array}{lr}
\dot{\hat{z}} = (I_r \otimes \hat{A} - \Gamma \otimes \hat{B}\hat{C}) \hat{z}  + (\mathcal{F} \otimes \hat{B}) u, \\
\hat{y} = 
(\mathcal{H} \otimes \hat{C}) \hat{z}.\\	
\end{array}
\right.
\end{equation}
where 
$$
\Gamma =  \begin{bmatrix}
\ihLambda_{1} & \\ & 0
\end{bmatrix}, \
\mathcal{F} = \begin{bmatrix}
\hat{F}_{1} \\  \tfrac{1}{\sqrt{n}}  \mathds{1}_n^\top F 
\end{bmatrix},\
\mathcal{H} = \begin{bmatrix}
\hat{H}_{1} & \tfrac{1}{\sqrt{n}} H  \mathds{1}_n 
\end{bmatrix}.
$$
Here, $\Gamma$ is not yet a Laplacian matrix,  but it has only one zero eigenvalue at the origin and all the other eigenvalues are positive real. To restore a network interpretation in the reduced-order model $\bm{\tilde{\Sigma}}$, the following theorem is provided in \cite{ChengAuto2019}, which states that there exits a similarity transformation between $\Gamma$ and an undirected graph Laplacian matrix.\\

\begin{thm} \label{thm:LapReal}
	A real square matrix $\Gamma$ is similar to the Laplacian matrix $\mathcal{L}$ associated with an weighted undirected connected graph, if and only if $\Gamma$ is diagonalizable and has an eigenvalue at $0$ with multiplicity $1$ while all the other eigenvalues are real and positive.
\end{thm}

By Theorem \ref{thm:LapReal}, we find a reduced Laplacian matrix $\hat{L}$ which has the same spectrum as $\Gamma$, namely, there exists a nonsingular matrix $\mathcal{T}_n$ such that
$
\hat{L} = \mathcal{T}_n^{-1}\Gamma\mathcal{T}_n.
$
The matrix $\hat{L}$ characterizes a reduced connected undirected graph $\hat{\mathcal{G}}$, which contains $r$ vertices.
%
%Now, denote $\Lambda = \text{diag} \{\lambda_1,\lambda_2,\cdots, \lambda_{r-1}, \lambda_{r} \}$ and consider the SVDs of $\mathcal{A}$ and $\hat{L}$:
%\begin{equation}
%	-\mathcal{A} = \mathcal{U}_a  \Lambda \mathcal{U}_a^\top, \ \hat{L} = \mathcal{U}_b  \Lambda \mathcal{U}_b^\top.
%\end{equation}
%with $\mathcal{U}_a^\top = \mathcal{U}_a^{-1}$ and $\mathcal{U}_b^\top = \mathcal{U}_b^{-1}$. 
%
%Then, it is easy to see that there is a nonsingular transformation matrix 
%\begin{equation}
%	\mathcal{T} =  \mathcal{U}_a \mathcal{U}_b^\top
%\end{equation}
%satisfying 
%$\hat{L} = -\mathcal{T}^{-1} \mathcal{A} \mathcal{T}$,
%where $ \mathcal{T}^{-1} = \mathcal{U}_b \mathcal{U}_a^\top  \mathcal{E}^{1/2}$. 
Applying a coordinate transform $\hat{z} = (\mathcal{T}_n \otimes I_{{r}} )\hat{x}$ to the system $\bm{\tilde{\Sigma}}$ in \eqref{sysRedMod} yields a reduced order network model
\begin{equation} \label{sysRedNet}
\bm{\hat{\Sigma}}: 
\left\{
\begin{array}{lr}
\dot{\hat{x}} = (I_r \otimes \hat{A} - \hat{L} \otimes \hat{B}\hat{C}) \hat{x}  + (\hat{F} \otimes \hat{B}) u, \\
\hat{y} = 
(\hat{H} \otimes \hat{C}) \hat{x},\\	
\end{array}
\right.
\end{equation}
with
$
\hat{F} = \mathcal{T}_n^{-1} \mathcal{F}  \ \text{and} \
\hat{H} = \mathcal{H} \mathcal{T}_n.
$
It can be verified that the reduced order network system $\bm{\hat{\Sigma}}$ in \eqref{sysRedNet} preserves   synchronization. Moreover, denote the transfer matrices of $\bm{\Sigma}$, $\bm{\hat{\Sigma}}$, ${\bf{\Sigma_s}}$, ${\bf{\hat{\Sigma}_s}}$, ${\bf{\Sigma_a}}$, and ${\bf{\hat{\Sigma}_a}}$ by $G$, $\hat{G}$, $T_s$, $\hat{T}_s$, $T_a$, and $\hat{T}_a$, respectively. 
The approximation error can be analyzed as follows. 
\begin{equation} \label{eq:ErrParts}
\begin{split}
\lVert G - \hat{G} \rVert_{\mathcal{H}_\infty} 
&= \lVert (T_s + T_a) - (\hat{T}_s  + \hat{T}_a) \rVert_{\mathcal{H}_\infty} \\
&\leq 
\lVert T_s - \hat{T}_s \rVert_{\mathcal{H}_\infty} +
\lVert T_a - \hat{T}_a \rVert_{\mathcal{H}_\infty},
\end{split}
\end{equation}
%The overall approximation error can be evaluated based on the reduction results of the stable system $\bf{\Sigma_s}$ and the average module $\bf{\Sigma_a}$.
in which \textit{a priori} upper bound on the reduction error of the stable system $\bf{\Sigma_s}$ is given as 
\begin{equation}\label{eq:errboundstab}
\lVert T_s - \hat{T}_s \rVert_{\mathcal{H}_\infty}  \leq  2 \sum\limits_{i=r}^{n-1}\sum\limits_{j=1}^{\ell}\sigma_i \tau_j + 2\sum\limits_{i=1}^{r-1}\sum\limits_{j=k+1}^{\ell}\sigma_i \tau_j,
\end{equation}
with $\sigma_i$ and $\tau_i$ the diagonal entries of $\Sigma_{\mathcal{G}}$ and  $\Sigma_{\mathcal{D}}$ in \eqref{eq:balancing}, respectively.
Denote $S_i$ and $\hat{S}_i$ as the transfer matrices of ${\bm{\Sigma}_i}$ and $\bm{\hat{\Sigma}}_i$, respectively. If $S_i - \hat{S}_i \in \mathcal{H}_\infty$, we obtain 
\begin{equation} \label{eq:errboundaverage}
\lVert T_a - \hat{T}_a \rVert_{\mathcal{H}_\infty} \leq \frac{1}{n} \lVert H\mathds{1}_n \mathds{1}_n^\top F \rVert_2 \cdot \lVert S_i - \hat{S}_i \rVert_{\mathcal{H}_\infty}.
\end{equation}

In several special cases, the \textit{a priori} error bounds on $\lVert G - \hat{G} \rVert_{\mathcal{H}_\infty}$ in \eqref{eq:ErrParts} can be obtained.
The first case is when we only reduce the dimension of the network while the agent dynamics are untouched as in \cite{Besselink2016Clustering}. In this case, we obtain $\lVert {T_a - \hat{T}_a} \rVert_{\mathcal{H}_\infty} = 0$, which yields 
\begin{equation}
\lVert G - \hat{G} \rVert_{\mathcal{H}_\infty} = 
\lVert T_s - \hat{T}_s \rVert_{\mathcal{H}_\infty} \leq 2 \sum_{i=r}^{n-1}\sum_{j=1}^{\ell} \sigma_i \tau_j.
\end{equation}
 
The second case is when  the average module is not observable from the outputs of the overall system $\bm{\Sigma}$ or not controllable by the external inputs. Consider 
\begin{equation} \label{eq:specialcase2}
H\mathds{1}_n= 0, \ \text{or} \ \mathds{1}_n^\top F = 0.
\end{equation}
Then, the approximation between $\bm{\Sigma}$ and $\bm{\hat{\Sigma}}$ is bounded by
$
\lVert G - \hat{G} \rVert_{\mathcal{H}_\infty}
= 	\lVert T_s - \hat{T}_s \rVert_{\mathcal{H}_\infty},
$
whose upper bound is given in \eqref{eq:errboundstab}. A special example of \eqref{eq:specialcase2}  can be found in \cite{Monshizadeh2013stability,Monshizadeh2014,leiter2017graph}, where the output matrix $H$ in \eqref{sysh} is taken as in \eqref{eq:outputAEP}. 
\vspace*{12pt}
 
\section{Conclusions}
\label{sec:conclusion}

In this chapter, model reduction techniques for linear dynamical networks with diffusive couplings have been reviewed. 
There exists a vast amount of literature on this topic, and the reference list in this chapter is certainly not complete. For example, in e.g., \cite{Rao2013graph,chu2011structured,Chow2013PowerReduction,Dorfler2013Kron,Monshizadeh2017Constant}, the approximation approaches are developed based on singular perturbation approximation and applied to reduce the complexity of chemical reaction networks and power networks. In \cite{jongsma2018CycleRemoval}, the interconnection topology is simplified by removing cycles in the network, and in e.g., \cite{biyik2008area,mlinaric2018synchronization,kawano2019data}, preliminary results for reducing nonlinear dynamical networks are developed.  Recently, a lot of interest is taken in the combination of network reduction and controller and observer designs. For example, \cite{xue2017lqg} presents a Linear Quadratic Gaussian (LQG) controller for large-scale dynamical networks using the clustering-based reduction, and \cite{sadamoto2017average,Umar2019average} proposes the average state observer based on reduced-order network models. 

Generally speaking, order reduction methods for linear network systems have been extensively investigated. However, the approximation of complex network containing nonlinear couplings or subsystems is still challenging, and  the existing results on nonlinear networks are far from satisfactory. Another challenge in this area is the order reduction of heterogeneous networks, i.e., network systems composed of nonidentical subsystems. 
 
\vspace*{20pt}

\bibliographystyle{abbrv} 
\bibliography{RedNet}
\end{document}